# Algorithmic Problems in the Braid Group

by

Elie Feder

A dissertation submitted to the Graduate Faculty in Mathematics in partial fulfillment of the requirements for the degree of Doctor of Philosophy, The City University of New York.

2003

This manuscript has been read and accepted by the Graduate Faculty in Mathematics in satisfaction of the dissertation requirement for the degree of Doctor of Philosophy.

\_\_\_\_\_\_\_\_\_\_\_\_\_ \_\_\_\_\_\_\_\_\_\_\_\_\_\_\_\_\_\_\_\_\_\_\_\_\_\_\_\_\_\_
Date       Chair of Examining Committee

\_\_\_\_\_\_\_\_\_\_\_\_\_ \_\_\_\_\_\_\_\_\_\_\_\_\_\_\_\_\_\_\_\_\_\_\_\_\_\_\_\_\_\_
Date       Executive Officer

<u>Michael Anshel</u>
<u>Al Vasquez</u>
<u>Burt Randol</u>

Supervisory Committee

The City University of New York

# Abstract

# ALGORITHMIC PROBLEMS IN THE BRAID GROUP

by

Elie Feder

Advisor: Professor Michael Anshel

The study of braid groups and their applications is a field which has attracted the interest of mathematicians and computer scientists alike. Their basic structure has been studied as far back as Gauss who considered the notion of a braid when studying the orbit of the first observed asteroid, Ceres. Besides for the value in studying the braid group in a theoretical framework, the braid group has been found to have diverse applications. While its applications to knot theory has been known for many years, its applicability to the field of cryptography has only been realized recently[AAG]. Since this time, a study of algorithmic problems in the braid group and their complexity have acquired a great practical significance, in addition to its intrinsic theoretical beauty.

We begin with a review of the notion of a braid group. We then discuss some known solutions to decision problems in braid groups. We then move on to proving new results in braid group algorithmics. We offer a quick solution to the generalized word problem in braid groups, in the special case of cyclic subgroups. We illustrate this solution and its complexity using a multitape Turing machine. We then turn to a discussion of decision problems in cyclic amalgamations of groups. Again using a multitape Turing machine, we solve the word problem for the cyclic amalgamation of two braid groups. We analyze its complexity as well.

We then turn to a more general study of the conjugacy problem in cyclic amalgamations. We revise and prove some theorems of Lipschutz[L1] and show their application to cyclic amalgamations of braid groups. We generalize this application to prove a new theorem regarding the conjugacy problem in cyclic amalgamations.

We then discuss some application of braid groups, culminating in a section devoted to the discussion of braid group cryptography. We conclude with a discussion of some open questions that we would like to pursue in future research.

# **<u>Acknowledgements</u>**

I would like to give a special thanks to my advisor, Michael Anshel, for all his guidance. I would like to thank Juan Gonzalez-Meneses for his help with some of the results in this paper. I would like to thank Joan Birman for her helpful discussions. I would like to thank Seymour Lipschutz and Bernie Domanski, whose works were helpful to my research.

# Table of Contents

# List of Figures

# **Outline and Summary of New Results**

For the reader who is familiar with the appropriate background material and wants an outline, and a survey of the *new* results found in this paper, the following few pages should suffice.

For a refresher on the assumed definitions and theorems, I would suggest two excellent references. For the necessary background on the braid group, see [B]. For the necessary background on combinatorial group theory, and more specifically on free products with amalgamation, see [MKS](especially chapter 4).

Section 1 of this paper is an introduction to the braid group, its presentations and its decision problems. It contains a survey of what is already known regarding these decision problems.

In section 2, we begin presenting new results. We consider a specific instance of the generalized word problem in the braid group; the case of cyclic subgroups. In other words, given two braid words $X, Y \in B_n$, we are asked to determine if there exists an integer c such that $X^c = Y$, i.e. if Y is in the cyclic subgroup generated by X. We have the following definition:

<u>Definition</u>: We define a function $\exp: B_n \to \mathbf{Z}$ where $\exp(\beta)$=the sum of the exponents of the generators which make up $\beta$.

We study this problem using a multitape Turing machine, following the lead of Anshel and Domanski [Do1,Do2,AD] who used such machines to study the complexity of certain decision problems in groups. We design and implement an algorithm, concluding section 2 with the following theorem:

<u>Theorem 1</u>: The generalized word problem in the case of a cyclic subgroup of the braid group $B_n$ is solvable in almost all cases. More specifically, if we are given

X,Y $\in B_n$, with exp(X)≠0, and if M is the maximal length of X and Y in BKL generators, then we can determine in $O(M^2n)$ whether Y is in the cyclic subgroup generated by X. We can also determine the specific power c such that $X^c=Y$.

In Section 3, we introduce the notion of a cyclic amalgamation and discuss its decision problems. We construct the cyclic amalgamation of two braid groups, a group which we use as a platform for new results in the rest of this paper.

In Section 4, we consider the word problem for the cyclic amalgamation of two braid groups. We show how its solution is dependent upon our Theorem 1; the solution to the generalized word problem for cyclic subgroups of the braid group. We design an algorithm to solve the word problem and implement it on a multitape Turing machine. After a complexity analysis, we conclude with the following theorem:

<u>Theorem 2</u>: Let G be the cyclic amalgamation of two braid groups, $B_1$ and $B_2$, each of braid index n. Furthermore, assume that the generators of the cyclic subgroups being amalgamated are powers of braid group generators,. Then G has a solvable word problem. More specifically, given a word w in G, where w=$b_1b_2…b_k$ such that alternating $b_i$ come from different factors, and where b is the maximal braid word length of the $b_i$'s, we can determine if w=1 in $O(b^2k^3n)$.

In section 5, we consider the conjugacy problem in cyclic amalgamations. We revisit some old results of Seymour Lipschutz[L1]. We revise some of his definitions and prove modified forms of some of his theorems. The main, new definitions and theorems are as follows:

<u>Definition</u>: We say that an element h in a group G is *conjugate-power-search-solvable* if for any w in G, we can decide whether or not w is conjugate to a power of h. Additionally, if w is, in fact, conjugate to a power of h, we can find which power of h it is conjugate to. If there are more then one, then there are finitely many and we can find them all.

Definition: We say that an element h in a group G is *double-coset-search-solvable* if for any pair u,v in G we can decide whether or not there exist integers r and s such that $h^r u h^s = v$.. Additionally if such an r and s exist, we can find them. If there are many such values, there are finitely many, and we can find them all.

Definition: We say that an element h in a group G is *super-semicritical* if h has the following properties:

 a′) h is non-self-conjugate

 b′) h is conjugate-power-search-solvable

 c′) h is double-coset-search-solvable.

Definition: We say that an element h in a group G is *super-critical* if in addition to properties a′), b′) and c′), we also have

 d′) If $h^m u = u h^m$, then u is a power of h.

Based upon these definitions we formulate and prove the following theorems:

Theorem 4′: Let A and B be groups with solvable conjugacy problems. Let a be a super-critical element of A and let b be a non-self-conjugate and power-search-solvable element of B. Then the free product of A and B amalgamating a and b has a solvable conjugacy problem.

Theorem 5′: Let G be the free product of two groups, A and B, with solvable conjugacy problem amalgamating a cyclic subgroup H generated by super-semicritical elements a and b in the factors. Then G has a solvable conjugacy problem.

 We then apply these theorems to the case of the cyclic amalgamation of two braid groups, and get the following nontrivial result.

<u>Theorem 6</u>: Let G be the free product of two braid groups amalgamating cyclic subgroups, each generated by a power of a generator in their respective braid group. Then G has solvable conjugacy problem.

Based upon the specific application of the theorems above to braid groups, we are led to consider a new generalization. We begin with some definitions.

<u>Definition</u>: We say that a group G is *exp-invarient*, if the sum of the exponents of its defining relations is equal to zero. (i.e. application of any of the relations does not effect the value of exp).

<u>Definition</u>: We say that an element h in G is *conjugate–search-solvable* if given u and v in G, we can decide if there exists an integer c such that $h^c u h^{-c} = v$, i.e. if u can be conjugated to v by a power of h. Additionally, we require that we can find all such values for c, if they exist.

We prove the following new theorem:

<u>Theorem 7</u>: Let A and B be exp-invarient groups with solvable conjugacy problem. Let $a \in A$ and $b \in B$ have the following properties:

i) $\exp(a) \neq 0$ and $\exp(b) \neq 0$

ii) a and b are conjugate-search-solvable.

Let G be the free product of A and B amalgamating the cyclic subgroups generated by a and b. Then G has a solvable conjugacy problem.

In section 6, we discuss applications of braid group algorithmics, with a specific discussion on its application to problems in knot theory.

Section 7 contains a thorough discussion of the latest field of applicability of braid group algorithmics—cryptography. From the time [AAG] showed how braid groups could be used in designing public key cryptosystems, the interest and research into all aspects of braid groups and their decision problems has greatly multiplied.

In section 8, we conclude this paper with a discussion of open problems and direction for further research into braid group algorithmics.

# 1 Introduction to Braids
## 1.1 The Braid Group

### 1.1.1 Definitions and Motivation for Study

Braids were studied as early as the beginning of the nineteenth century by Gauss[Ga]. He considered the notion of a braid when studying the orbit of the first observed asteroid, Ceres [E]. Since then, the notion of a braid and the braid group have evolved significantly. We will begin by discussing the underlying geometric motivation for the formation of the braid group and then lead into the formal algebraic definition of this group in terms of its presentations.

To study the braid group geometrically, we begin with some definitions and examples.

<u>Definition</u>: A braid is a set of n disjoint strings in 3-space, all attached to a horizontal bar at the top and at the bottom. Each string always heads downwards as we move from the top bar to the bottom bar. (In general we will not draw the top and bottom bar.)

<u>Example</u>:

: 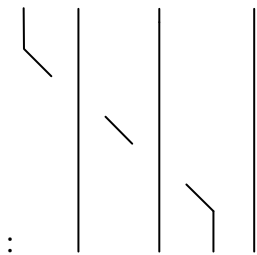

Braid of 4 strings

<u>Definition</u>: Two braids are considered equivalent if we can rearrange the strings in the two braids to look the same without passing any strings through one another or themselves.

Note: We must keep the bars fixed and the strings attached to the bars.

Example:

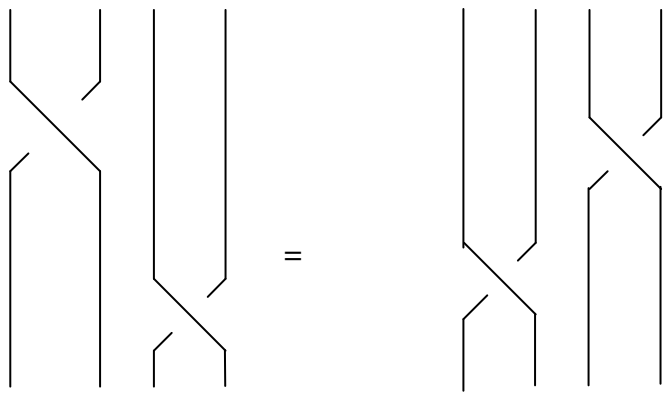

Two equivalent braids of four strands each

## 1.1.2 Formation of the braid group geometrically

Consider the set of all braids with a fixed number of strings. Multiplication of two braids with the same braid index (number of strings) is defined by concatenation, i.e. place one braid under the other.

Example:

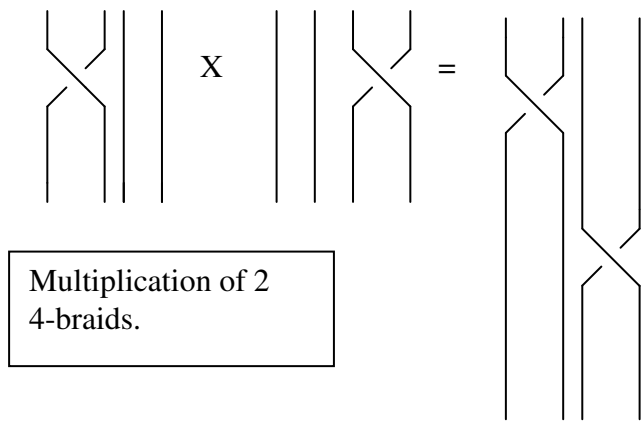

Multiplication of 2 4-braids.

The identity braid is simply the braid consisting of strings which do not cross.

Example: 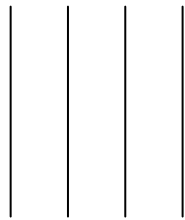

Identity braid on 4 strings

We can always find the inverse of a braid by "untwisting" all twists. The inverse of a braid can be attained by reflecting the braid with respect to a horizontal line.

Example:

Inverse of 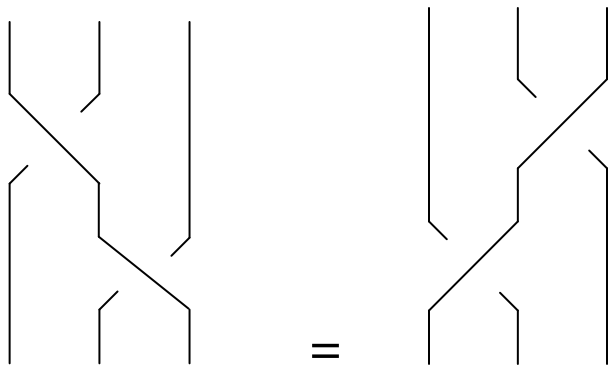

It can be verified that this multiplication is associative and, thus, the set of n-strand braids forms a group, termed the braid group of index n. In order to study certain properties of the braid group and to find algorithmic solutions to its decision problems, we must deal with the braid group in a more concrete sense. We will, therefore, consider some presentations for the braid group.

## 1.2 Presentations of the Braid Group

<u>Definition</u> A presentation of a group is a way of specifying a group completely in terms of a set of generators and a set of defining relations on these generators.

## 1.2.1 Artin Presentation of the Braid Group

The first presentation for the braid group was put forth by Emil Artin in 1925[A]. He defined $B_n$, the braid group with braid index n, by the following generators and relations:

<u>Generators</u>: $\sigma_1, \sigma_2, \ldots, \sigma_{n-1}$

<u>Defining relations</u>: (1) $\sigma_s \sigma_t = \sigma_t \sigma_s$, for $|t-s| > 1$

(2) $\sigma_s \sigma_t \sigma_s = \sigma_t \sigma_s \sigma_t$, for $|t-s| = 1$

$\sigma_i$ represents the braid in which the $(i+1)^{st}$ string crosses over the $i^{th}$ string while all other strings remain uncrossed.

<u>Example</u>:

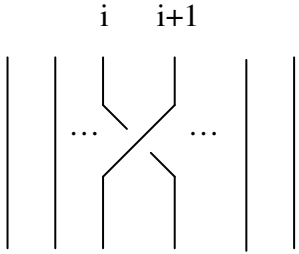

The above diagram represents $\sigma_i$

Example of relation (1):

$$\sigma_1 \sigma_3 = \sigma_3 \sigma_1$$

=

Example of relation (2):

$$\sigma_1 \sigma_2 \sigma_1 = \sigma_2 \sigma_1 \sigma_2$$

=

These generators and relations completely describe $B_n$.

## 1.2.2 BKL-Presentation of the Braid Group

An alternate presentation was suggested in [BKL]. Instead of considering the generators as braids in which two adjacent strings cross, they considered the generators as braids in which any two strings cross. Namely, the generators are: $a_{ts}$, $1 \leq s < t \leq n$, where $a_{ts}$ represents the braid in which the $t^{th}$ string crosses over the $s^{th}$ string while the $s^{th}$ and $t^{th}$ string cross in front of all intermediate strings.

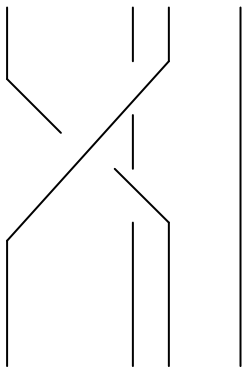

Example: $a_{31} \in B_4$

Notice that $\sigma_i = a_{(i+1,i)}$. Thus, the Artin generators are a subset of the BKL generators.

Now let us consider the BKL defining relators:

1) $\qquad a_{ts}a_{rq} \qquad = \qquad a_{rq}a_{ts} \qquad$ if $(t-r)(t-q)(s-r)(s-q) > 0$

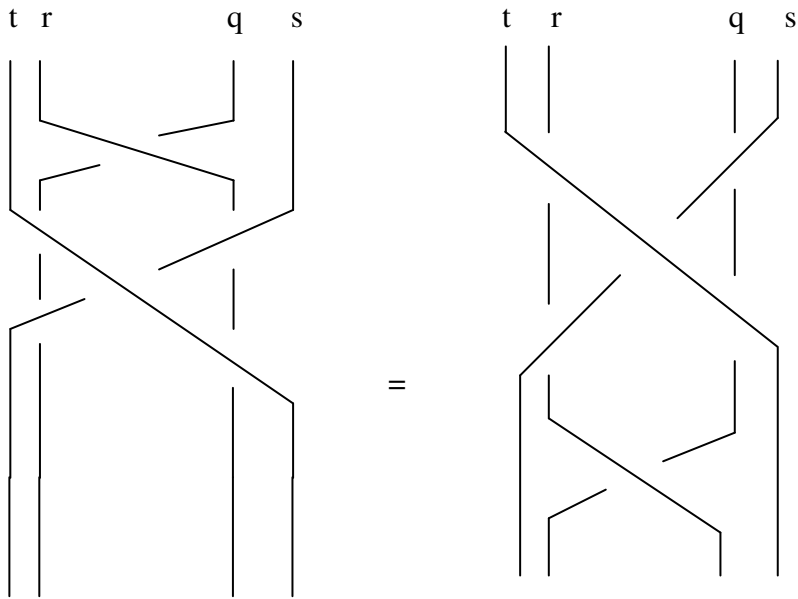

This means that $a_{ts}$ and $a_{rq}$ commute if t and s do not separate r and q.

2)    $a_{ts}a_{sr}$ = $a_{tr}a_{ts}$ = $a_{sr}a_{tr}$

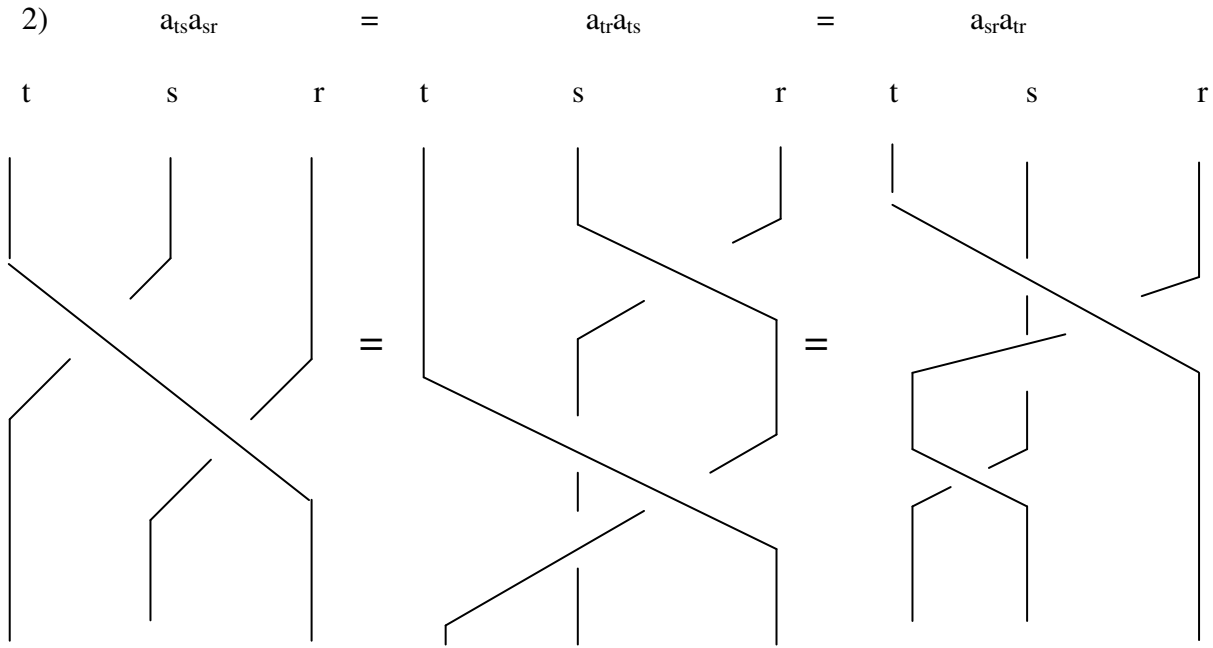

This is a partial commutativity when the two generators share a strand.

Note: We have not drawn in the intermediate strings, but they may exist. If they do, they are all crossed over by t, s and r

## 1.3 Decision Problems

Whenever we have a group presentation, we are led to considering "decision problems" for the group. This notion was introduced by Max Dehn in 1911[D1]. Some of these problems include:

### 1.3.1 The Word Problem

Given a group G defined by a presentation, and given a word W in G, given as a product of the generators and their inverses, can we decide in a finite number of steps whether or not W defines the identity of G? In other words, can we decide if W can be reduced to the identity in G by applying the defining relations in G?

Note: This problem is equivalent to the problem of deciding whether or not two words in G are equivalent under the defining relations. For, to say that $W_1=W_2$ is equivalent to saying that $W_1W_2^{-1} = 1$.

There have been a number of approaches towards solving the word problem for the braid group.

i) Artin's Solution

The first one to offer a solution to the word problem was the person who posed the problem—Emil Artin, in [A]. He considered the map $\varphi: B_n \to S_n$, where $S_n$ represents the symmetric group of order n. The map sends $\sigma_i$ to the transposition (i, i+1). He then considered the kernel of this map, which is called $P_n$, the pure braid group on n strings. By studying $P_n$, he showed how to put any braid into a normal form called a *combed braid*. He thus attained a solution to the word problem in the braid group. It seems that his solution is exponential in the word length[BKL].

ii) Garside-Thurston Solution

In 1969, Garside[G] suggested a different solution to the word problem in the braid group. This solution is exponential in braid index and word length. In 1992, Thurston[T] improved upon Garside's approach and proved that for word length $|W|$, there exists an algorithmic solution to the word problem in $B_n$ with complexity $O(|W|^2 n \log n)$. We will now outline their basic approach. We begin with some definitions.

<u>Definition</u>: A *positive braid* in $B_n$ is any braid which is expressed as a product of positive powers of the generators $\sigma_1,...,\sigma_{n-1}$. We can form the monoid of all positive braids; we call this monoid $B_n^+$.

We then introduce the braid word $\Delta = (\sigma_1 \sigma_2 ... \sigma_{n-1})(\sigma_1 \sigma_2 ... \sigma_{n-2})...(\sigma_1 \sigma_2)(\sigma_1)$. This is called the *fundamental braid*. It can be shown that $\Delta^2$ generates the center of $B_n$. Garside proved that any braid could be represented by a word of the form $\Delta^r P$, where r is the maximal integer for all such forms and where P is a positive braid word. In such a normal form, we say that r=inf of the given braid. Garside had no unique representation for P. This is where Thurston improved Garside's algorithm. He showed that P can be uniquely factored as a product $P_1 P_2 ... P_k$, where each $P_i$ is a specific positive braid called a *permutation braid*. They have thus produced a normal form for any braid word. This solves the word problem. For, given two words W and V, to determine if W=V, we simply must convert W and V to their respective normal forms and then check if they are equivalent.

iii) The Birman-Ko-Lee Solution

In 2001, [BKL] improved upon the algorithm of Garside-Thurston by introducing the [BKL] presentation for the braid group which is discussed above. Using these new generators and defining relations, they mirrored the algorithm of Garside-Thurston in this new context. They defined a new fundamental word

$\delta = a_{n(n-1)}a_{(n-1)(n-2)}\ldots a_{21}$

( $= \sigma_{n-1}\sigma_{n-2}\ldots\sigma_2\sigma_1$ in the Artin generators). They then showed that any word $W \in B_n$ can be represented uniquely in the BKL generators by a word of the form $\delta^j A_1 A_2 \ldots A_k$ where $A = A_1 A_2 \ldots A_k$ is a positive word, j is maximal and k is minimal for all such representations. The $A_i$'s are positive braids which are uniquely determined by their associated permutations. They call these $A_i$'s *canonical factors*. The above form of A is shown to have what is called the left-weighted property. The increase in the efficiency of this algorithm over that of Thurston-Garside is that there are much fewer canonical factors then there are permutation braids. The number of different canonical factors is the $n^{th}$ Catalan number, which is less then $4^n$. On the other hand, the number of permutation braids is n!, which grows faster then $k^n$ for any k. The [BKL] algorithm for the word problem has complexity $O(|W|^2 n)$

iv) Dehornoy's Solution

In [De], Dehornoy introduces a completely different approach to solving the word problem. This algorithm is based upon the existence of a linear ordering on braids and a method of reducing braids, *handle reduction*. Although it is hard to make an exact comparison between this method and those mentioned above, Dehornoy claims that in practice this algorithm is more efficient than all previously known methods.

### 1.3.2 The Generalized Word Problem

Given a group presentation for G, and given H, a subgroup of G, can we decide in a finite number of steps whether or not a given word W in G is also contained in H?

Note: The word problem is the special case of the generalized word problem when H is equal to the trivial subgroup in G.

In [BD], it is shown that $B_5$ contains an isomorphic copy of the direct product of two free groups of rank 2, as a subgroup. Additionally, in [Mih], it was shown that the direct product of two free groups has unsolvable generalized word problem. This, being the case, one can conclude that the generalized word problem for braid groups of index ≥ 5, is unsolvable. However, in section 2 we will prove that for the specific case of cyclic subgroups, whose generator has nonzero exponential sum, the generalized word problem is, in fact, solvable. We will design an algorithm and implement it on a multitape Turing machine, and thereby analyze the complexity of this problem.

### 1.3.3 The Conjugacy Problem

Given a group presentation for G, and given two group words U and V, can we decide in a finite number of steps whether or not U and V are conjugate in G? In other words, does there exist an element W in G such that $W^{-1}UW = V$ ?

The conjugacy problem is a significantly more difficult problem then the word problem, in the case of braid groups. A number of solutions to the conjugacy problem have been suggested over the years.

i) Garside's solution

In [G], Garside proved that the conjugacy problem for braid groups is solvable. For a given braid P, he defined a set called the *summit set* of P. This set consists of all conjugates of P which have maximal inf. (Recall, if the normal form of P is $\Delta^r P_1 P_2 \ldots P_k$, then we say that inf(P)=r and sup(P)=r+k.). Clearly, the summit set is a conjugacy invariant. Thus, to determine if two braids are conjugate, it suffices to compare their summit sets. The difficulty lies in efficiently computing the summit set of a given element.

ii) Elrifai and Morton's solution

Garside's solution was improved by Elrifai and Morton in [EM]. They defined a set known as the *super summit set*(SSS). Given a braid word W, the SSS of W is the set of elements which are conjugate to W that have both minimal sup and maximal inf. The SSS is clearly a subset of the summit set, and is easier to compute.

iii) BKL's solution

Subsequently, [BKL] improved this algorithm even further by studying the SSS under their presentation of $B_n$, as opposed to Artin's.(see Section 1.2) Despite these improvements, all known solutions are exponential in the length of the words involved. However, [BKL] conjectured that the solution is polynomial in word length of the elements being compared. This conjecture has not been proven as of yet. However, in [FG1], the authors improved the algorithm and made some computations which give numerical support for the conjecture.

iv) linear algebraic solutions

Another possible approach to solving the conjugacy problem in braid groups is a linear algebraic approach. Instead of trying to solve the conjugacy problem directly in the braid group, we solve it in a *representation* of the braid group. A representation of a group is a homomorphism from the group to $GL_k(R)$, the set of k x k invertible matrices over R, for some integer k and some R, a domain with characteristic zero. The idea is as follows: to determine if two braid elements are conjugate, we (i) find their image under a given representation. (ii) determine if these images are conjugate in $GL_k(R)$ (iii) lift the solution back to the braid group.

In order to carry out a linear algebraic solution to the conjugacy problem, we need to find "good" representations for $B_n$. One characteristic that we would like a representation to have is that it is faithful (i.e.1-1). If a faithful representation of $B_n$ exists, then we say that $B_n$ is linear. We want our representation to be faithful so that when we lift the solution from $GL_k(R)$ back to the braid group, we are assured of getting one braid, the one desired. A long time candidate for a faithful representation of the braid group was the Burau representation[Bu]. But, in [Mo] it was proven that in fact it was not faithful for $n \geq 9$. Later, it was shown to be unfaithful for $n \geq 5$.

In 1999, Krammer[K1] constructed a representation which he proved to be faithful for n=4. This was later proven to be faithful for all n,[Bi],[K2] and therefore, $B_n$ is now known to be linear for all n. This representation is known as the Lawrence-Krammer representation, being that Krammer's construction was based upon earlier work by Lawrence[L]. There have been many attempts at solving the conjugacy problem in the braid group using, both, the Burau representation and the Lawrence–Krammer representation. We will hold off from a discussion of these approaches until Section 7, when we discuss the application of the conjugacy problem in the braid group to cryptography.

### 1.3.4   The Shortest Word Problem

Given a group presentation G and a word W in G, consider the class of all words in G which are equivalent to W, under the defining relations of G. Can we find the element in this class of words which has the shortest word length in the generators and their inverses? In 1991, [PR] showed that the shortest word problem in the braid group is at least as hard as an NP-complete problem.

## 2 Generalized Word Problem for Cyclic Subgroups in the Braid Group

### 2.1 The Problem

We consider the generalized word problem for the class of cyclic subgroups of $B_n$. Given a braid x and a braid y, both in $B_n$, can we find an algorithm to determine whether y is in the cyclic subgroup generated by x or not? If it is, can we find the exponent k such that $x^k=y$ and demonstrate this equivalence? What is the running time for this algorithm?

Note: This problem was raised in the newsgroup sci.crypt[S] and attracted much interest.

### 2.2 Some Lemmas

To solve this problem, we introduce a definition and some lemmas.

Definition: We define a function $\exp:B_n \rightarrow \mathbf{Z}$ where $\exp(\beta)$=the sum of the exponents of the generators which make up $\beta$.

Example: $\exp(\sigma_1 \sigma_3^{-3} \sigma_2^2 \sigma_1)= 1-3+2+1 = 1$.

Lemma 1: exp is invariant under braid equivalence. In other words if u=v in $B_n$, then exp(u)=exp(v).

Proof: Let u=v in $B_n$. Then we can get from u to v by a chain of elements in $B_n$, where each element can be attained from the previous element by applying one of the defining relations in $B_n$. Thus, it suffices to show that applying any of the defining relations in $B_n$ will not change the value of exp. Let us, therefore, consider the defining relations in $B_n$.

1. $\sigma_i\sigma_j=\sigma_j\sigma_i$ where $|i-j|>1$. Well, $\exp(\sigma_i\sigma_j) = 2 = \exp(\sigma_j\sigma_i)$.

2. $\sigma_i\sigma_{i+1}\sigma_i=\sigma_{i+1}\sigma_i\sigma_{i+1}$. Well, $\exp(\sigma_i\sigma_{i+1}\sigma_i) = 3 = \exp(\sigma_{i+1}\sigma_i\sigma_{i+1})$.

<u>Note</u>: If we have an occurrence of $\sigma_i\sigma_i^{-1}$ or $\sigma_i^{-1}\sigma_i$ and we replace it with 1, we also do not alter the value of exp. This is because $\exp(\sigma_i\sigma_i^{-1}) = \exp(1) = \exp(\sigma_i^{-1}\sigma_i) = 0$. The lemma is thus proven.

<u>Lemma 2</u>: exp is additive. Namely $\exp(u\cdot v)=\exp(u)+ \exp(v)$

<u>Proof</u>: Since $u\cdot v$ is simply the concatenation of the words u and v, the sum of the exponents of the new word $u\cdot v$ is simply the sum of the exponents of u plus the sum of the exponents of v. Any cancellation caused by the concatenation is simply the removal of $\sigma_i\sigma_i^{-1}$ or $\sigma_i^{-1}\sigma_i$, both of whose exponent sum is zero. Thus removal of these terms will not effect the total exponent sum.

<u>Lemma 3</u>: $\exp(\beta^{-1})= -\exp(\beta)$

<u>Proof</u>: Let $\beta=\sigma_{i1}^{p1}\sigma_{i2}^{p2}\ldots\sigma_{in}^{pn}$. Then $\exp(\beta)= p_1+p_2+\ldots+p_n$. Now $\beta^{-1}=\sigma_{in}^{-pn}\ldots\sigma_{i2}^{-p2}\sigma_{i1}^{-p1}$. Then $\exp(\beta^{-1})=-p_1-p_2-\ldots-p_n=-(p_1+p_2+\ldots+p_n)= -\exp(\beta)$.

<u>Lemma 4:</u> exp is invariant under conjugacy in the braid group. In other words, if u is conjugate to v in $B_n$, then $\exp(u)=\exp(v)$.

<u>Proof</u>: Let us assume that $u = \beta^{-1}v\beta$, where $\beta\in B_n$.

By lemmas 2 and 3, we have $\exp(u) = \exp(\beta^{-1}v\beta)= \exp(\beta) + \exp(v) + \exp(\beta^{-1})= \exp(v)$.

<u>Lemma 5</u>: If $y=x^k$, then $\exp(y)= k\cdot\exp(x)$.

<u>Proof</u>: This follows directly from the additivity of exp (lemma 2).

## 2.3 The Solution

a) Now if we are given that $y=x^k$ for some k, and if we further assume that $\exp(x)\neq 0$, then we can easily solve for k. Namely, $k= \exp(y)/\exp(x)$.

b) Now, let us assume that we are given x, y in $B_n$, with $\exp(x) \neq 0$ and we are told to determine if $y = x^k$ for some k. In other words, is y in the cyclic subgroup generated by x? We proceed as follows. Assume it is. Then $y = x^k$, for some k. We can compute the only possible value for k. Namely $k = \exp(y)/\exp(x)$. Now compute $x^k$. Decide if $y = x^k$ by the known algorithms for the word problem[CKLHC].

Note: The above solution is limited to the case where $\exp(x) \neq 0$. For any braid word x, where $\exp(x) = 0$, a different solution would be required. As $\exp(x) = 0$ for any $x \in [B_n, B_n]$, the commutator of $B_n$, such a solution would certainly be desirable.

## 2.4 The Outline of the Algorithm

Given: $x, y \in B_n$.

Question: Is y in the cyclic subgroup generated by x?

Step 1: Compute $\exp(x)$ and $\exp(y)$

Step 2: Compute $\exp(y)/\exp(x) = k$

Step 3: Compute $x^k$

Step 4: Decide $x^k = y$? If yes, then y is in the cyclic subgroup generated by x and we know the power, k. If no, then y is not in the cyclic subgroup generated by x.

Note: In the case that $\exp(y) = 0$, then $k = 0$. By convention, define $x^0 = 1$. Thus, this reduces to the word problem, as a subcase.

## 2.5 Programming Language for Multitape Turing Machine

We will follow the model of some papers by Domanski and Anshel[Do1,Do2,AD] for our programming language for multitape-Turing machines. We assume a familiarity with the basic model of a multitape-Turing machine which consists of an input tape, k worktapes, and an output tape. Let the symbol being read by tape head A on tape X be denoted as X(A). Moving tape head A on X one tape cell to the right will be denoted by A=A+1. Similarly, moving tape head A on X one tape cell to the left will be denoted by  A=A-1. We allow standard high-level programming constructs such as IF-THEN-ELSE, DO-WHILE and DO-UNTIL… One can verify that these constructions can be converted to standard Turing machine notation. We also allow the number theoretic operations mod and div.  In the course of the algorithm, we will denote the beginning of an explanatory comment with /* , and the end of such a comment with  */.

## 2.6  The Formal Algorithm GWP(X,Y)

<u>Input</u>: X,Y$\in$ B$_n$ , exp(X)$\neq$0.

<u>Task</u>: Determine if there exists an integer c such that X$^c$=Y.

<u>Output</u>: GWP(X,Y)= "X$^c$=Y", if Y is the c-th power of X

"Y is not a power of X", if there does not exist an integer c such that X$^c$=Y.

We will use a multitape-Turing machine with three worktapes; X,Y and Z.

Tape X will contain the braid word X=X$_1$X$_2$…X$_n$  followed by the $ symbol.

Tape Y will contain the braid word Y=Y$_1$Y$_2$…Y$_m$  followed by the $ symbol.

X$_1$,…,X$_n$,Y$_1$,…,Y$_m$ each denote a generator of the braid group or its inverse.

We will define sign(X$_i$)=1 if X$_i$ is a generator and sign(X$_i$)=-1 if X$_i$ is the inverse of a generator.

Tape Z will contain the integers.

A will be a tape head situated on tape X and begins under the first letter of X, namely, X$_1$ . B will be a tape head situated on tape Y and begins under the first letter of Y, namely Y$_1$.

EX and EY will be tape heads situated on tape Z and will begin under zero.

/* EX and EY will be used to keep track of the exponent sum of X and Y */

Do Until X(A)=$
    If sign(X(A))=1,
        Then Do;

```
                    EX=EX+1
              End;
        Else Do;
              EX=EX-1
            End;
        A=A+1
End; /* of DO-UNTIL  */
```

/* Z(EX) is now equal to exp(X), the exponent sum of X.  */

```
Do Until Y(B)=$
      If sign(Y(B))=1,
            Then Do;
                    EY=EY+1
              End;
        Else Do;
              EY=EY-1
            End;
        B=B+1
End  /*  of DO-UNTIL  */
```

/* Z(EY) is now equal to exp(Y), the exponent sum of Y.  */

If Z(EY) mod (Z(EX)) =0,  /* exp(X) divides exp(Y)  */
      Then Do
              c= Z(EY) div (Z(EX))
              Form U=$X^c$

/* U is formed by concatenating c copies of X. This can be done in linear time on a

multitape Turing machine as is done in [Do2] */

/* Comparison(G,H) is assumed to be a routine which, given G,H$\in B_n$, determines if G=H

under the relations of $B_n$. The output is True if G=H and False if G$\neq$H   See [CKLHC]. If

we assume that L is the minimal of the canonical lengths of G and H, then this routine

runs in $O(L^2 n \log n)$ if G and H were given in Artin generators and in $O(L^2 n)$ if G and H

were given in BKL generators. */

              If Comparison(U,Y)=True
                    Then Do
                          Output "$X^c$=Y".
                    End;

```
            Else Do;
                    Output "X is not a power of Y"
            End;
Else Do;
        Output "X is not a power of Y"
    End;
End;
```

## 2.7  Complexity Analysis of GWP(X,Y)

We will let M be the maximal length of X and Y. The computation of both Z(EX) and Z(EY) is bounded by M. Forming $X^c$ is linear in the length of $X^c$. What is the length of $X^c$? Well, c is bounded by M and, therefore, the length of $X^c$ is certainly bounded by $M^2$. As mentioned above, if X and Y are given in BKL generators, then the routine Comparison(G,H) is bound by $O(L^2n)$, where L is the minimal length of G and H. In our case where we are comparing $X^c$ and Y, this minimal length is bound by M because the length of Y is bound by M. Thus, Comparison($X^c$,Y) is bound by $O(M^2n)$. Thus, the entire algorithm GWP(X,Y) is bound by $O(M^2n)$.

We can now conclude with the following theorem:

<u>Theorem 1</u>: The generalized word problem in the case of a cyclic subgroup of the braid group $B_n$ is solvable in almost all cases. More specifically, if we are given

$X, Y \in B_n$, with exp(X)≠0, and if M is the maximal length of X and Y in BKL generators, then we can determine in $O(M^2n)$ whether Y is in the cyclic subgroup generated by X. We can also determine the specific power c such that $X^c=Y$.

## 2.8 The Root Problem in the Braid Group

A problem which is very closely related to the generalized word problem in the case of cyclic subgroups of the braid group, is the root problem in the braid group. In a sense, these are opposite problems. The root problem is as follows: Given a braid word $w \in B_n$, and a positive integer r, we are asked to determine if w has any $r^{th}$ roots. In other words, does there exist a braid $b \in B_n$, such that $b^r=w$?

In [St], it was shown that this problem is solvable. However, no complexity bound was found for this solution. In [Si], this solution was generalized to small Gaussian groups, a generalization of braid groups which was introduced in [DP]. It would be nice if we could find a bound on the complexity of these algorithms.

## 3  Cyclic Amalgamations

### 3.1  Definitions and Decision Problems

Let A and B be groups. Let $h \in A$ and $k \in B$. Assume that <h>, the cyclic subgroup generated by h is isomorphic to <k>, the cyclic subgroup generated by k, under the isomorphism $\varphi$. Let C=<h>=$\varphi$(<k>). Let G=A$*_C$B denote the free product of A and B with <h> and <k> amalgamated by identifying h with k. We call G a *cyclic amalgamation* of A and B. In terms of generator and relations, if we let G(A) be the generators of A, G(B) the generators of B, R(A) the relations in A, R(B), the relations in B, then G=A $*_C$ B = <G(A),G(B); R(A),R(B), h=k >. When we discuss cyclic amalgamations in the course of this paper, we will understand the definition in terms of generator and relations. Cyclic amalgamations are but one instance of a general construction of free products with amalgamation.

When studying cyclic amalgamations, as in the case with any group, we can ask ourselves a few questions:

1  Can we solve the word problem? What is its complexity?

2  Can we solve the conjugacy problem? What is its complexity?

3  Can we solve the normal form problem? What is its complexity?

These problems have all been considered in different cases. We will show later that the word problem is solvable iff the generalized word problem for the relevant cyclic subgroups in each factor are solvable. The normal form problem is always solvable in cyclic amalgamations, as in the more general case of free products with amalgamation. Namely, for G=A $*_C$ B , any element $g \in G$ can be represented uniquely as g=ca$_1$b$_1$a$_2$b$_2$…a$_n$b$_n$, where c∈C, each a$_i$≠1 ∈ A, and each b$_i$≠1 ∈ B (see [Mi2] for instance). The conjugacy problem, on the other hand, is a problem whose solution is not so simple. In [L3] it was shown that the conjugacy problem for the cyclic amalgamation of free groups is solvable. However, in [Mi1], an example was given of a free product of two free groups amalgamating *finitely generated* subgroups which has an unsolvable conjugacy problem. In [CT], it was shown, that if A and B are sixth-groups, and h∈A and k∈B have the same order, then the free product of A and B amalgamating the cyclic subgroups generated by h and k has solvable conjugacy problem. In [L1] and [L2], the conjugacy problem for cyclic amalgamations was discussed at length and many results were given, some of which were not proven. Later in this paper we will discuss some of these results and prove them. We will also prove other theorems on this topic. We will show that braid groups are an interesting setting to discuss cyclic amalgamations in. The case which we will consider is as follows:

Let $B_1$ and $B_2$ be braid groups on n strings.

i.e. $B_1 = <\sigma_1,…,\sigma_{n-1}; \sigma_i\sigma_{i+1}\sigma_i=\sigma_{i+1}\sigma_i\sigma_{i+1}, \sigma_k\sigma_j=\sigma_j\sigma_k \text{ for } |k-j|>1 >$

and $B_2 = <\tau_1,…,\tau_{n-1}; \tau_i\tau_{i+1}\tau_i=\tau_{i+1}\tau_i\tau_{i+1}, \tau_k\tau_j=\tau_j\tau_k \text{ for } |k-j|>1 >$

Let G be the cyclic amalgamation of $B_1$ and $B_2$.

i.e. $G = B_1*_C B_2 = <\sigma_1,…,\sigma_{n-1},\tau_1,…,\tau_{n-1}; R(B_1), R(B_2), \sigma_k^p=\tau_j^r >$

for any k,j={1,…,n-1} and p,r positive integers.

## 4 Algorithm for the Word Problem in the Cyclic Amalgamation of Two Braid Groups

### 4.1 The Solution

<u>Problem</u>. Given a word w in G, can we determine if w=1?

<u>Solution</u>: Assume that $w=b_1b_2...b_k$ where each $b_i$ and $b_{i+1}$ come from a different factors.($B_1$ and $B_2$ are referred to as *factors* in G since G is the free *product* of $B_1$ and $B_2$ with amalgamation.) We can assume this, for if $b_i$ and $b_{i+1}$ are both in $B_i$, then we can "collapse" them to get $b_ib_{i+1}$ in $B_i$. We thus have a shorter word in G. We can continue this process until the above condition is met. Now assume, without loss of generality, that $b_1$ is in $B_1$. Using the algorithm for the generalized word problem for cyclic subgroups of braid groups( we will call this algorithm "GWP" ), check if $b_1$ is in the cyclic subgroup of $B_1$ generated by $\sigma_k^p$. Namely, apply the routine GWP($\sigma_k^p$,$b_1$). If not, check the analogous condition for $b_2$ with regard to $B_2$. If none of the $b_i$'s are in the associated cyclic subgroups, then we can conclude that $w \neq 1$, by a well-known lemma. (see [Mi1], for instance) If $b_i$ is in $<\sigma_k^p>$, for instance, then let us assume that $b_i=\sigma_k^{pc}$. Then replace $b_i$ with $\tau_j^{rc}$, by the added defining relation in G. We thus have $b_{i-1}\tau_j^{rc}b_{i+1} \in B_2$. We have thus reduced the word length of w. We now continue this process by induction until we reach 1 or until we cannot reduce any further. If we reach 1 then we conclude w=1. If we get to a point where we can not reduce any further, then we conclude that $w \neq 1$.

## 4.2 The Formal Multitape Turing Machine Algorithm, I(w)

<u>Input</u>: Let $w= b_1b_2...b_n$ where all $b_i$ are nontrivial braid words in $B_1$ or $B_2$ and $b_i$ and $b_{i+1}$ are from different factors.

<u>Task</u>: To determine if w=1 in G.

<u>Output</u>: " I(w)=True " if w=1

" I(w)= False " if $w \neq 1$

We will use the same programming language for multitape Turing machines that was discussed in section 2.5. However, instead of the cells of the worktapes containing braid generators and their inverses, the cells in this setup will contain braid words, from $B_1$ or $B_2$. If we were to be more rigorous, we would let the cells contain braid generators and attain the braid words from this setup. However, our assumption will not effect the overall complexity of the algorithm.

Consider two worktapes W and U, with heads H and R on tape W and head L on tape U.

/*   W will be the tape which we will be reducing w on, while U will be where we save our results- the reduced form of w.   */

The initial configuration for the Turing machine will be as follows:

```
Tape W:   b₁ b₂…bₙ1
          H  R
Tape U:   1
          L
```

where H is positioned under the leftmost symbol of W

      R is positioned under the second from the left symbol of W

      L is positioned under the symbol 1 on U

<u>Note</u>: 1 refers to the identity in the braid group.

/* We will move along W and test if we can reduce the length of W by using the relation provided by the amalgamated subgroup. We do this by testing if each of the $b_i$'s are contained in the amalgamated subgroup. Since this subgroup is cyclic (with the generator having nonzero exponent sum), we can test this using the routine GWP, as above.

We assume that the routine GWP(y,x) takes x and y as its input, and determines if $y=x^c$ for some integer c. If so, then the routine outputs c. We assume that $x^0=1$, the identity braid word.   */

```
Do Until W(H)=1
    If W(H)∈ B₁
       Then Do;
            If  GWP(W(H),σₖᵖ)=c   /* if W(H)=σₖᵖᶜ */
                Then Do; /* Replace by the element τⱼʳᶜ in B₂ and reduce. */
                     W(H)=U(L)τⱼʳᶜW(R)
                     If W(R)≠1
                         Then R=R+1
```

```
                          If U(L)≠1
                              Then L=L-1
                       End;
                  Else Do;   /* Move further along W  */
                      L=L+1
                      U(L)=W(H)
                      Do Until H=R;
                          H=H+1
                      End;
                      If W(R)≠1
                          Then R=R+1
                  End;
     Else Do      /* if W(H)∈ B₂ */
         If GWP(W(H),τⱼʳ)=d     /* if W(H)=τⱼʳᵈ  */
     Then Do; /*  Replace by the element σₖᵖᵈ in B₁ and reduce.*/
                         W(H)=U(L)σₖᵖᵈW(R)
                         If W(R)≠1
                             Then R=R+1
                         If U(L)≠1
                             Then L=L-1
                    End;
         ElseDo;   /*  Move further along W  */
                      L=L+1
                      U(L)=W(H)
                      Do Until H=R;
                          H=H+1
                      End;
                      If W(R)≠1
                          Then R=R+1
         End; /* of ELSE-DO */
     End;  /* of ELSE-DO  */
```

/* If in the course of our reduction, we end up with a word in the amalgamated subgroup, then this algorithm will throw us into an infinite loop, converting this word from an element of $B_1$ to an element of $B_2$. to an element of $B_1$… We, therefore, must add an If-Then clause to terminate the program in such a case and thereby get us out of such a loop */

```
        If ( W(R)=1 and U(L)=1)
             Then Do
                  If W(H)=1
                      Then output " I(W)=True"     /* W=1  */
                  Else output "I(W)=False"        /*  W≠1  */
```

```
                    Terminate program.
                End;
        End  /* of If-Then  */
End;  /* of Do-Until  */
If  U(L)=1,
        Then output "I(W)=True"  /*  W=1.  */
Else output "I(W)=False"   /*  W≠1  */
End;
```

## 4.3  Complexity Analysis

What is the complexity of this algorithm? The most significant step in this algorithm is the application of the routine GWP. We apply this many times. As proven in Theorem 1 in Section 2.7, the complexity for GWP is $O(M^2 n)$ where M is equal to max {length of generator of the cyclic subgroup, length of word being checked}. Since any power of $\sigma_k^p$ or $\tau_j^r$ has a greater length then these cyclic generators, if the length of the word being checked is less then that of the generator, we can conclude that this word is not such a power without running GWP. Thus, we can assume that the length of the word being checked is greater then that of the generator. How long are the words being checked? Well, with each "collapse" the new word being checked is longer then before. What's the longest possible word length? Well, let us assume that the maximum (braid group) length of the elements $b_1, b_2, \ldots, b_k$ is b. Now when, we replace a word $b_i$ by a power of the generator, we can only decrease its length. (Why? Firstly, any replacement must keep the exp constant. Now a positive power of a generator has the shortest word length of all words of equivalent exp, being that all its terms are of the same sign.) Thus the length of any word we must check will certainly be <bk . Thus, each application of the GWP algorithm has complexity bound $O(b^2 k^2 n)$. How many times must we apply this algorithm? Notice that the amount of applications of GWP is equal to the amount of values that W(H) takes on in the running of this algorithm. How many such values are there? Well, we begin with tape W having k possible values for H. Each time GWP yields false, H=H+1, and we "use up" one of these k values. Thus we have at most k GWP applications which yield false, before the algorithm will terminate with W(H)=1. How many applications of GWP can yield true before the algorithm terminates? Well, each time GWP yields true, we reduce w and the length decreases by at least 1. Thus, the bound on the number of these reductions is k, the starting length of w. Thus we have a total bound of 2k applications of GWP. We can conclude as follows:

Theorem 2: Let G be the cyclic amalgamation of two braid groups, $B_1$ and $B_2$, each of braid index n. Furthermore, assume that the generators of the cyclic subgroups being amalgamated are powers of braid group generators,. Then G has a solvable word problem. More specifically, given a word w in G, where w=$b_1b_2...b_k$ such that alternating $b_i$ come from different factors, and where b is the maximal braid word length of the $b_i$'s, we can determine if w=1 in $O(b^2k^3n)$.

## 4.4 Generalization

In proving that the word problem in G is solvable, the only property of braid groups we used is the following:

i)  Given an element b∈$B_n$, exp(b)≠0, we can determine if a word w in $B_n$ is a power of b and if so, we can find the power.

Whenever property i) holds for an element b in a group G, we say that b is *power-solvable*. If, additionally, we can determine the power, we say that b is *power-search-solvable*.

We can thus generalize to a well-known theorem.(see [L4], for instance)

Theorem 3: Let G be the free product of groups amalgamating cyclic subgroups which are generated by power-search-solvable elements in the factors. Then G has a solvable word problem.

# 5  Conjugacy Problem in Cyclic Amalgamations

## 5.1  Lipschutz's Theorems

We now turn towards considering the conjugacy problem in G. In treating the conjugacy problem in the cyclic amalgamation of two braid groups, I would like to show that the solution is merely one example of a more general theorem regarding the solution of the conjugacy problem in cyclic amalgamations. I believe that this general theorem was already discussed years ago, but was never made precise enough to apply to braid groups. In [L1], Lipschutz discusses properties of groups and elements within them which are necessary to ensure that their cyclic amalgamation will have a solvable conjugacy problem. In order to discuss his results, we must begin with some definitions from his paper.

<u>Definition</u>: We say that an element h in a group G is *non-self-conjugate* if its distinct powers are in different conjugacy classes.

<u>Definition</u>: We say that an element h in a group G is *power-solvable* if for any w in G, we can decide whether or not w is a power of h. (as above)

<u>Definition</u>: We say that an element h in a group G is *conjugate-power-solvable* if for any w in G, we can decide whether or nor w is conjugate to a power of h.

<u>Definition</u>: We say that an element h in a group G is *double-coset-solvable* if for any pair u,v in G we can decide whether or not there exist integers r and s such that $h^r u h^s = v$.

<u>Definition</u>: We say that an element h in a group G is semicritical if h has the following properties:

    a) h is non-self-conjugate

    b) h is conjugate-power-solvable

    c) h is double-coset-solvable.

<u>Definition</u>: We say that an element h in a group G is *critical* if in addition to properties a), b) and c), we also have

    d) If $h^m u = u h^m$, then u is a power of h.

We can now state two of Lipschutz's theorems for which he stated no proof.

Theorem 4: Let A and B be groups with solvable conjugacy problems. Let u be a critical element of A and let v be a non-self-conjugate and power-solvable element of B. Then the free product of A and B amalgamating u and v has a solvable conjugacy problem.

Theorem 5: Let G be the free product of groups with solvable conjugacy problem amalgamating a cyclic subgroup generated by semicritical elements in the factors. Then G has a solvable conjugacy problem.

In the mathematical review[A], the reviewer wrote "The reviewer awaits with interest the proof of theorem 3(which in our paper is called theorem 5) and some significant nontrivial new examples to give the properties discussed in this paper some substance." In [Hu], the author showed that the standard embedding of a free product with amalgamation into the corresponding HNN extension is conjugacy preserving. He then showed that the in the HNN extension which corresponds with the cyclic amalgamation in Lipschutz's theorem 3, the conjugacy problem in solvable. He, thereby, proved theorem 3. Additionally, some examples were provided by Lipschutz in [L2]. However, I believe that theorem 2 (in our paper, theorem 4) has not been demonstrated, nor has a direct proof for theorem 3 appeared.

## 5.2 Modification and Proof of Lipschutz's Theorems

I would like to give proofs for these two theorems of Lipschutz. However, I believe that some of his conditions must be strengthened. Additionally, I believe that the case of braid groups provides a significant nontrivial new example for theorem 5, which the reviewer was looking for to give the theorem some substance.

Firstly, let's add some definitions.

Definition: We say that an element h in a group G is *conjugate-power-search-solvable* if for any w in G, we can decide whether or not w is conjugate to a power of h. Additionally, if w is, in fact, conjugate to a power of h, we can find which power of h it is conjugate to. If there are more then one, then there are finitely many and we can find them all.

Definition: We say that an element h in a group G is *double-coset-search-solvable* if for any pair u,v in G we can decide whether or not there exist integers r and s such that $h^r u h^s = v$. Additionally if such an r and s exist, we can find them. If there are many such values, there are finitely many, and we can find them all.

Definition: We say that an element h in a group G is *super-semicritical* if h has the following properties:

a′) h is non-self-conjugate

b′) h is conjugate-power-search-solvable

c′) h is double-coset-search-solvable.

Definition: We say that an element h in a group G is *super-critical* if in addition to properties a′), b′) and c′), we also have

d′) If $h^m u = u h^m$, then u is a power of h.

We now state modified versions of theorems 4 and 5 which we will prove. We will also show that theorem 5′ has the cyclic amalgamation of braid groups as a nontrivial example.

Theorem 4′: Let A and B be groups with solvable conjugacy problems. Let a be a super-critical element of A and let b be a non-self-conjugate and power-search-solvable element of B. Then the free product of A and B amalgamating a and b has a solvable conjugacy problem.

Theorem 5′: Let G be the free product of two groups, A and B, with solvable conjugacy problem amalgamating a cyclic subgroup H generated by super-semicritical elements a and b in the factors. Then G has a solvable conjugacy problem.

We will first prove theorem 5′ and then prove theorem 4′ in a similar manner.

Proof of theorem 5′: We need to prove the following: Given u, v cyclically reduced elements of G. We can decide if u is conjugate to v in G.

The algorithm will be based on Solitar Theorem [MKS- theorem 4.6] which reads as follows:

Theorem: Let $G = A *_H B$. Then every element of G is conjugate to a cyclically reduced element of G. Moreover, suppose that g is a cyclically reduced element of G. Then:

(i) If g is conjugate to an element h in H, then g is in some factor and there is a sequence $h, h_1, h_2, \ldots, h_t, g$ where each $h_i$ is in H and consecutive terms of the sequence are conjugate in a factor.

(ii) If g is conjugate to an element g′ in some factor, but not in a conjugate of H, then g and g′ are in the same factor and are conjugate in that factor.

(iii) If g is conjugate to an element $p_1 \ldots p_r$ where $r > 1$ and $p_i$, $p_{i+1}$ as well as $p_1$, $p_r$ are in distinct factors, then g can be obtained by cyclically permuting $p_1 \ldots p_r$ and then conjugating by an element of H.

By the above theorem, we can assume u and v are cyclically reduced. If not, conjugate until we get a cyclically reduced element. The proof then breaks up into 3 cases.

Case 1: u and v have different lengths in G. Then they are not conjugate.

Case 2: u and v both have length 1. By Solitar's Theorem, either

(i) u and v are in the same factor and are conjugate in that factor. If so the problem reduces to the conjugacy problem in A and B whose solutions are assumed to be known. Or

(ii) u and v are in different factors.(assume without loss of generality that $u \in A$ and $v \in B$) Then there exists a sequence of elements $u, h_1, h_2, \ldots, h_s, v$ such that u is conjugate to $h_1$ in A, $h_1$ is conjugate to $h_2$ in B,…, $h_s$ is conjugate to v in B, where $h_1, \ldots, h_s$ are in H. But different elements in H are not conjugate to each other because a and b are non-self-conjugate. Thus, we have the sequence u, $h_I$, v , where $h_I$ is in H.

Now given u, we can check if it is conjugate to an element in H because a and b are assumed to be conjugate-power-search-solvable. If u is conjugate to a power of a, find this power, say m. Now test if v is conjugate to $a^m$ (by the conjugacy problem in B, which is assumed to be solvable). If it is, then u is conjugate to v in G. If not, not.

<u>Case 3</u>: u, v have length n>1. Assume $u=u_1u_2\ldots u_n$ and $v=v_1v_2\ldots v_n$. By Solitar's Theorem, u is conjugate to v iff there exists an i and j such that

$h^m u_i\ldots u_n u_1\ldots u_{i-1} h^{-m} = v_j\ldots v_n v_1\ldots v_{j-1}$ i.e. iff a permutation of u is conjugate to a permutation of v by an element in H. Now, we have n permutations of u and n permutations of v. We thus have $n^2$ cases to test. We must, therefore, find a method to test if

(*) $h^m a_1\ldots a_n h^{-m} = b_1\ldots b_n$ for any $a_1\ldots a_n$ and $b_1\ldots b_n \in G$

Proceed as follows: Decide if h and $a_1$ commute. This is possible by the solution to the word problem in A and B.

<u>Case (a)</u>: Assume not. Then, $b_1^{-1} h^m a_1 = h^n$ for some n, by properties of amalgamated products. ( Namely, multiply (*) by $b_1^{-1}$ on both sides. The length in G will only remain equivalent on both sides if the above condition is met) This implies that

$h^m a_1 h^{-n} = b_1$. Now since h is double-coset-search-solvable, we can find if there exists such an m. Find all values of m and n which satisfy the above equation. Then test if

$h^m a_1\ldots a_n h^{-m} = b_1\ldots b_n$. for each such m.

<u>Note</u>: We can test this because the word problem in G is solvable. This is true because h is power-solvable (which follows from the fact that h is conjugate-power-solvable) and we can, therefore, apply theorem 3(from Section 4.3) to it.

If the above equation holds for some such m, and for some $a_1\ldots a_n$ and $b_1\ldots b_n$ permutations of u and v, then u is conjugate to v. If not, not.

Case (b): Assume that $\sigma_k$ and $a_1$ commute. Then $b_1^{-1}a_1=h^c$ for some c. This implies that $a_1=b_1h^c$. We can solve for c since h is power-search-solvable. Now substitute into the equation $h^m a_1\ldots a_n h^{-m}=b_1\ldots b_n$ and use the fact that h and $a_1$ commute to get

$b_1 h^m h^c a_2\ldots a_n h^{-m} = b_1 b_2\ldots b_n. \to h^m(h^c a_2)a_3\ldots a_n h^{-m} = b_2\ldots b_n$. This is the same as above but we have a shorter word length in G. Thus, we can check this by induction. This proves theorem 5′.

Proof of theorem 4′: Cases 1 and 2 are identical to theorem 5′.

Case 3: u, v have length n>1. Assume $u=u_1u_2\ldots u_n$ and $v=v_1v_2\ldots v_n$. By Solitar's Theorem, u is conjugate to v iff there exists an i and j such that $h^m u_i\ldots u_n u_1\ldots u_{i-1} h^{-m} = v_j\ldots v_n v_1\ldots v_{j-1}$ i.e. iff a permutation of u is conjugate to a permutation of v by an element in H. Now, we have n permutations of u and n permutations of v. We thus have $n^2$ cases to test. We must, therefore, find a method to test if $h^m a_1\ldots a_n h^{-m}=b_1\ldots b_n$, for any $a_1\ldots a_n$ and $b_1\ldots b_n$, permutations of u and v. Proceed as follows: We can assume, without loss of generality, that $a_1\in A$. [Why? If it weren't, then since $a_1\ldots a_n$ is cyclically reduced, $a_n\in A$. If so, we could decide if $u^{-1}$ is conjugate to $v^{-1}$(which is equivalent to deciding if u is conjugate to v), in which case we would be deciding if $h^m a_n^{-1}\ldots a_1^{-1} h^{-m}=b_n^{-1}\ldots b_1^{-1}$, in which case $a_n^{-1}\in A$(this is an observation of Lipschutz in [L2])]. Now decide if $a_1$ commutes with a (the generator of H in A) This is possible by solution to word problem in A. If it does, then since a is super-critical, we have that $a_1$ is a power of a, by property (d′). This is a contradiction to the fact that u is reduced, with length n. Thus, we can conclude that $a_1$ does not commute with h. Then, $b_1^{-1}h^m a_1=h^n$ for some n, by properties of amalgamated products. This implies that $h^m a_1 h^{-n}= b_1$. Now since a is double-coset-search-solvable, we can find if there exists such an m and n. Find all values which satisfy the above equation. Then test if $h^m a_1\ldots a_n h^{-m}=b_1\ldots b_n$ for each such m. If the above equation holds for some such m, and for some permutations of u and v, then u is conjugate to v. If not, not. Theorem 4′ is thus proven.

## 5.3 Application to the Braid Group

Theorem 6: Let G be the free product of two braid groups amalgamating cyclic subgroups each generated by a power of a generator in their respective braid group. Then G has solvable conjugacy problem.

Proof: We will show that we could apply theorem 5′ to this case. We already know that the conjugacy problem in the braid group is solvable. Thus, we are left to showing that in any braid group, a power of a generator is super-semicritical, i.e., that it satisfies properties (a′), (b′) and (c′).

Let $\sigma_k^p$ be the generator of the cyclic subgroup of the braid group which we are amalgamating. Let's take the properties one at a time:

(a′) $\sigma_k^p$ is non-self-conjugate.

Proof: We must show that distinct powers of $\sigma_k^p$ are in different conjugacy classes. Well, we showed above that elements in the same conjugacy class have the same value of exp. But, all powers of $\sigma_k^p$ have a different exp. Therefore, they are in different conjugacy classes.

(b′) $\sigma_k^p$ is conjugate-power-search-solvable.

Proof: Given a w in $B_n$, we must decide whether or not w is conjugate to a power of $\sigma_k^p$. If it is, we must find that power ( there can be at most one because $\sigma_k^p$ is non-self-conjugate). Well, assume that w is conjugate to $\sigma_k^{pc}$ for some integer c. Then $\exp(w)=\exp(\sigma_k^{pc})=pc$. $\to$ $c=\exp(w)/p$, which we can compute. We, therefore, must decide if w is conjugate to $\sigma_k^{\exp(w)}$, which can be decided using the solution to the conjugacy problem in braid groups.

(c′) $\sigma_k^p$ is double-coset-search-solvable.

Proof: Given a pair a and b in $B_n$, we must decide whether or not there exist integers m and n such that $\sigma_k^{pm} a \sigma_k^{pn} = b$.

Recall: In the proof of theorem 3′, this property was only necessary in the case where a and h (in this case a and $\sigma_k^p$) did not commute. Thus, we can assume in this proof that a and $\sigma_k^p$ do not commute, i.e. a is not in the centralizer of $\sigma_k^p$.

Well, if such an m and n do exist, then $\sigma_k^{pm} a \sigma_k^{-pm} \sigma_k^{p(m+n)} = b$

$\rightarrow \sigma_k^{pm} a \sigma_k^{-pm} = b \sigma_k^{-p(m+n)} \rightarrow \exp(a) = \exp(b) - p(m+n) \rightarrow m+n = (\exp(b) - \exp(a))/p$

Let us denote $b\sigma_k^{-p(m+n)}$ as d. Then $\sigma_k^{pm} a \sigma_k^{-pm} = d$. Thus our original question is equivalent to deciding if a is conjugate to d by a power of $\sigma_k^p$. We thus consider a more general question. Can we can decide if two elements a and b in $B_n$ are conjugate by a power of a given generator?

Lemma: Given a,b$\in B_n$, and a generator $\sigma_i$, there is a finite algorithm to determine whether there exists a k such that $\sigma_i^{-k} a \sigma_i^k = b$. Additionally, we can compute such a k, if it exists.

Proof: In [EM], the following is proven: Let a,b$\in B_n$. Let p=min{inf(a),inf(b)}. If $u^{-1}au=b$ and the canonical form for u is $U_1…U_n$, then $\inf((U_1…U_I)^{-1}a(U_1…U_i)) \geq p$ for all i$\leq$n. Applying this result to our case, we can conclude that if $\sigma_I^{-k} a \sigma_i^k = b$, for k>0, then $\inf(\sigma_i^{-j} a \sigma_i^j) \geq p$ for all j$\leq$k. Now in [G] it is proven that the number of words of fixed exp and inf$\geq$p is finite. Since all conjugates have the same exp(by lemma 4 in section 2.2), there are only a finite number of elements conjugate to a which have inf$\geq$p. Thus, we can proceed as follows: Compute $\sigma_I^{-1} a \sigma_I$, then $\sigma_I^{-2} a \sigma_I^2$, then $\sigma_I^{-3} a \sigma_I^3$,…until either (1) we get an element whose inf is less then p or (2) we enter a loop or (3) we get b. In case (3) the answer to the problem is positive and we have found the value for k. In case (1) and (2), we have determined that a is not conjugate to b by a positive power of $\sigma_i$. We still must test negative powers. Namely, perhaps $\sigma_i^k a \sigma_I^{-k} = b$, for k>0. If so then $\sigma_i^{-k} b \sigma_i^k = a$. We can then follow the same process as above with the places of a and b reversed. As noted above, this is a finite process. We have thus produced a finite algorithm to test if a is conjugate to b by a power of a given generator. Additionally, we can find the desired power using this algorithm. The lemma is proven.

We now get back to our problem. We left off with the task of determining if there exists an m such that $\sigma_k^{pm} a \sigma_k^{-pm} = d$ where a and d are elements in $B_n$, Thus apply the lemma to test if a is conjugate to d by a power of $\sigma_k$. If not, then there is no such m and thus there do not exist an m and n such that $\sigma_k^{pm} a \sigma_k^{pn} \neq b$ and we're done. If yes, find the value j such that $\sigma_k^{-j} a \sigma_k^{j} = d$. Is j a multiple of p? If not, again there is no such m and n and we're done. If yes, then we have found a value of m such that $\sigma_k^{pm} a \sigma_k^{-pm} = d$. Namely m=j/p. Therefore, there exist an m and n such that $\sigma_k^{pm} a \sigma_k^{pn} = b$. Namely, m=j/p and n=(m+(exp(a)-exp(b)))/p, as above. Now, we still have one consideration. To prove that $\sigma_k^p$ is double-coset-search-solvable, we must find all values of m and n such that… Perhaps, there is a q>j such that $\sigma_k^{-q} a \sigma_k^{q} = d$. Perhaps this will produce another possible m and n which we must test. How do we know we can stop at j? Well, assume this is true. Then $\sigma_k^{-j} a \sigma_k^{j} = \sigma_k^{-q} a \sigma_k^{q} \rightarrow \sigma_k^{q-j} a \sigma_k^{j-q} = a. \rightarrow a^{-1} \sigma_k^{q-j} a = \sigma_k^{q-j}$. This implies that a is in the centralizer of $\sigma_k^{q-j}$ where q-j>0. But by results in [FG2] it can be shown (see lemma below) that $\sigma_k$ and $\sigma_k^c$ have the same centralizer for any value of c. Thus, since a is in the centralizer of $\sigma_k^{q-j}$, a is in the centralizer of $\sigma_k$. This implies that a is in the centralizer of $\sigma_k^p$. This contradicts our assumption above. Therefore, we can conclude that there is no such q. Therefore, we have found the only possible values of m and n. We have, thereby, shown that $\sigma_k^p$ is double-coset-search-solvable.( whenever $\sigma_k^p$ and a don't commute)

We have, thus, shown that $\sigma_k^p$ is super-semicritical. We are now justified in applying theorem 5′ to the case of braid groups and the theorem is thus proven.

<u>Lemma</u>: For any generator $\sigma_k \in B_n$ and any c>0, the centralizer of $\sigma_k$ is the same as the centralizer of $\sigma_k^c$.

Proof: Well, in [FG2], it was shown that the generators of the centralizer of an element are completely determined by the edges of its minimal conjugacy graph. It, therefore, suffices to show that the edges of the minimal conjugacy graph for $\sigma_k$ are the same as those for $\sigma_k^c$. Now, by definition, the vertices of the minimal conjugacy graph of an element a are those positive elements (in $B_n$) which are conjugate to a. The edges of the graph are labeled by minimal simple elements (see the paper for definitions) as follows: An edge labeled by s goes from a vertex u to a vertex v if s is a minimal simple element and $s^{-1}us=v$. Here's an example from the paper.[FG2]

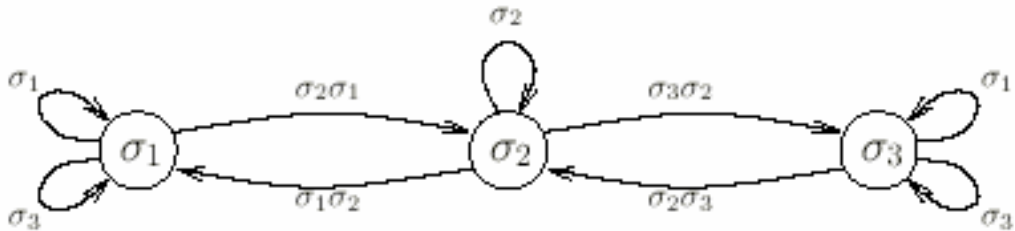

Figure 3: Minimal conjugacy graph of $\sigma_1 \in B_4^+$.

This example could be generalized for any $\sigma_i \in B_n$, for any n, as follows: The positive conjugates of $\sigma_i$ are $\sigma_1, \sigma_2, \ldots \sigma_{n-1}$. These are, therefore, the vertices of the minimal conjugacy graph for any $\sigma_i$. How are the edges of the graph labeled? Well, consider any vertex $\sigma_j$. Notice that for any k such that $|j-k|>1$ or j=k, we have $\sigma_k^{-1}\sigma_j\sigma_k=\sigma_j$. Therefore, we have an edges from $\sigma_j$ to itself labeled by all such $\sigma_k$. Now for $|j-k|=1$, we have $(\sigma_k\sigma_j)^{-1}\sigma_j(\sigma_k\sigma_j) = \sigma_j^{-1}\sigma_k^{-1}\sigma_j\sigma_k\sigma_j = \sigma_j^{-1}\sigma_k^{-1}\sigma_k\sigma_j\sigma_k = \sigma_k$. We thus have an edge labeled $\sigma_k\sigma_j$ going from $\sigma_j$ to $\sigma_k$ for all such k. This is how the minimal conjugacy graph for $\sigma_i$ looks, and the centralizer for $\sigma_i$ is determined by its edges. It is trivial to see that all the edges are labeled by minimal simple elements with respect to $\sigma_i$. Now let us consider the minimal conjugacy graph for $\sigma_i^c$. The positive conjugates of $\sigma_i^c$ are

$\sigma_1^c, \sigma_2^c, \ldots, \sigma_{n-1}^c$. How are the edges of the graph labeled? Well, consider any vertex $\sigma_j^c$. Notice that for any k such that $|j-k|>1$ or $j=k$, we have $\sigma_k^{-1}\sigma_j^c\sigma_k=\sigma_j^c$. Therefore, we have edges from $\sigma_j^c$ to itself labeled by all such $\sigma_k$. Now for $|j-k|=1$, we have $(\sigma_k\sigma_j)^{-1}\sigma_j^c(\sigma_k\sigma_j) = \sigma_j^{-1}\sigma_k^{-1}\sigma_j^c\sigma_k\sigma_j$

<u>Claim</u>: $\sigma_j^{-1}\sigma_k^{-1}\sigma_j^c\sigma_k\sigma_j = \sigma_k^c$.

<u>Proof</u>: By induction on c. For c=1, we showed this above. Assume true for c-1.

i.e. $\sigma_j^{-1}\sigma_k^{-1}\sigma_j^{c-1}\sigma_k\sigma_j = \sigma_k^{c-1}$. Now consider $\sigma_j^{-1}\sigma_k^{-1}\sigma_j^c\sigma_k\sigma_j = \sigma_j^{-1}\sigma_k^{-1}\sigma_j^{c-1}\sigma_j\sigma_k\sigma_j$

$= \sigma_j^{-1}\sigma_k^{-1}\sigma_j^{c-1}\sigma_k\sigma_j\sigma_k = \sigma_k^{c-1}\cdot\sigma_k = \sigma_k^c$ and the claim is proven. We thus have an edge from $\sigma_j^c$ to $\sigma_k^c$ which is labeled $\sigma_k\sigma_j$. It is trivial to see that all the edges in this minimal conjugacy graph are also labeled by minimal simple elements with respect to $\sigma_i^c$. Now notice that the edges in this minimal conjugacy graph for $\sigma_i^c$ are the same as those in the minimal conjugacy graph for $\sigma_i$. We can thus conclude that the centralizer of $\sigma_I$ is equal to the centralizer of $\sigma_I^c$ for any c>0. For c<0, the result follows from noticing that the centralizer for $a^c$ is equal to the centralizer for $a^{-c}$. The lemma is proven. This lemma is proven by a different method in [FRZ].

## 5.4  A New Generalization

Now that we found that the case of braid groups provides an interesting example for theorem 3′, we will show that we can generalize this example to state another theorem regarding the solvability of the conjugacy problem for cyclic amalgamations. We begin with some definitions.

<u>Definition</u>: We say that a group G is *exp-invariant*, if the sum of the exponents of its defining relations is equal to zero. (i.e. application of any of the relations does not effect the value of exp).

As shown, in lemma 3 in section 2.2, the braid group is exp-invariant. This property of braid groups was essential in proving that theorem 3′ applied to braid groups.

Note: There are many other examples of groups with this property. See [G], for instance.

Definition: We say that an element h in G is *conjugate–search-solvable* if given u and v in G, we can decide if there exists an integer c such that $h^c u h^{-c} = v$, i.e. if u can be conjugated to v by a power of h. Additionally, we require that we can find all such values for c, if they exist.

Note: The condition that h is conjugate-search-solvable is weaker then the condition of theorems 4′ and 5′ that h is double-coset-search-solvable.

As shown above, any power of any generator in a braid group is conjugate-search-solvable.

Theorem 7: Let A and B be exp-invariant groups with solvable conjugacy problem. Let $a \in A$ and $b \in B$ have the following properties:

i) $\exp(a) \neq 0$ and $\exp(b) \neq 0$

ii) a and b are conjugate-search-solvable.

Let G be the free product of A and B amalgamating the cyclic subgroups generated by a and b. Then G has a solvable conjugacy problem.

Proof: We will prove that a and b are super-semicritical and we can, therefore, apply theorem 5′. We will show that a has properties (a′),(b′) and (c′). (b follows similarly)

(a′) a is non-self-conjugate.

Proof: Since A is exp-invariant, and conjugating does not change the value of exp, we can conclude that two conjugate elements in A have the same value of exp. Now, for any i,

$\exp(a^i) = i \cdot \exp(a)$. Since $\exp(a) \neq 0$, for $i \neq j$, $\exp(a^i) \neq \exp(a^j)$. Thus, different powers of a are not conjugate to one another.

(b′) a is conjugate-power-search-solvable

<u>Proof</u>: Given w in A, we must decide whether or not w is conjugate to a power of a. If it is, we must find that power ( there can be at most one because a is non-self-conjugate). Well, assume that w is conjugate to $a^c$ for some integer c. Then since A is exp-invariant, $\exp(w)=\exp(a^c)=c\cdot\exp(a)$. $\rightarrow c=\exp(w)/\exp(a)$, which is nonzero and computable. We, therefore, must decide if w is conjugate to $a^c$. This can be decided using the solution to the conjugacy problem in A.

(c′) a is double-coset-search-solvable.

<u>Proof</u>: Given a pair u and v in A, we must decide whether or not there exist integers m and n such that $a^m u a^n = v$. Well, if there are, then $a^m u a^{-m} a^{(m+n)} = v \rightarrow a^m u a^{-m} = v a^{-(m+n)}$

$\rightarrow \exp(u) = \exp(v) - (m+n)\cdot\exp(a) \rightarrow m+n = (\exp(v)-\exp(u))/\exp(a)$

Let us denote $v a^{-(m+n)}$ as d. Then $a^m u a^{-m} = d$. Thus our original question is equivalent to deciding if u is conjugate to d by a power of a. This is decidable since a is conjugate-search-solvable. Therefore a and b are super-semicritical. We can now apply theorem 5′ to prove theorem 7.

# 6  Applications of Braid Groups

When studying any group, we can ask all of the above decision problems. But in the case of the braid group, they have unique applications.

<u>6.1  The Knot Equivalence Problem and The Unknotting Problem</u>

Descriptively, a knot is a piece of string in 3-space, which is allowed to arbitrarily wrap around itself until it eventually connects up to where it began, thus forming a closed curve. The Knot Equivalence Problem asks us to determine algorithmically when a given two knots are equivalent. Descriptively, we say two knots are equivalent if one can be deformed to the other in 3-space without cutting either one. This notion is formalized by the Reidemeister moves[R]. Namely, we say that knots K and J are equivalent if we can form a sequence of knots

K= $K_1 \to K_2 \to \ldots \to K_i \to \ldots \to K_n$ =J, where each $K_{i+1}$ can be attained from $K_i$ by a Reidemeister move. There are three types of Reidemeister moves as illustrated below:

Reidemeister move type 1

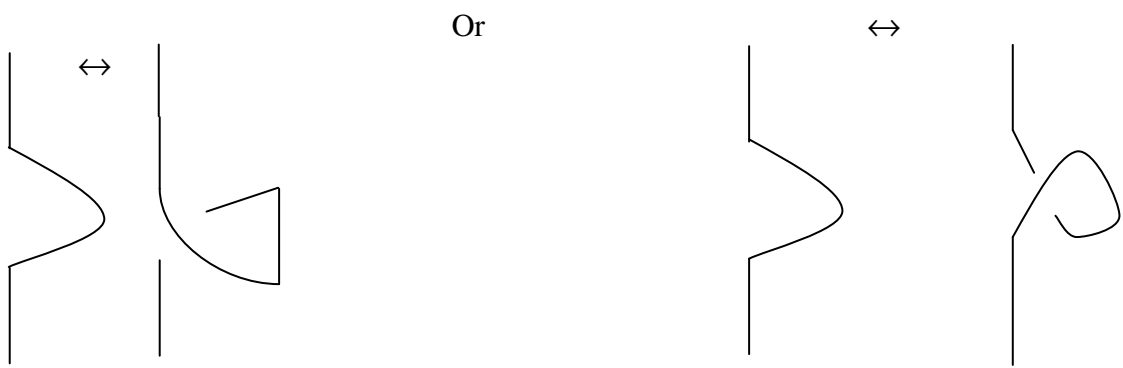

Reidemeister move type 2

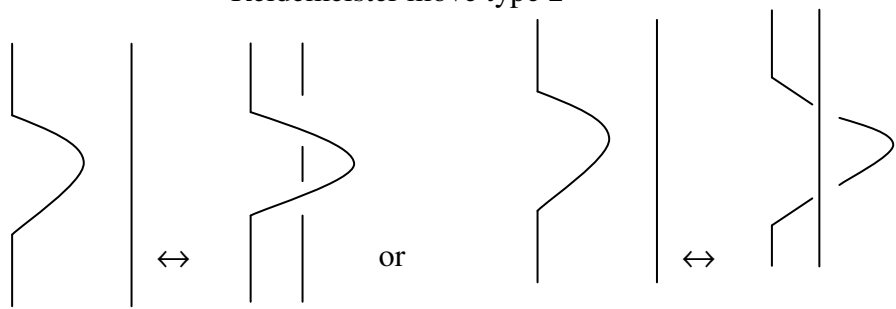

Reidemeister move type 3

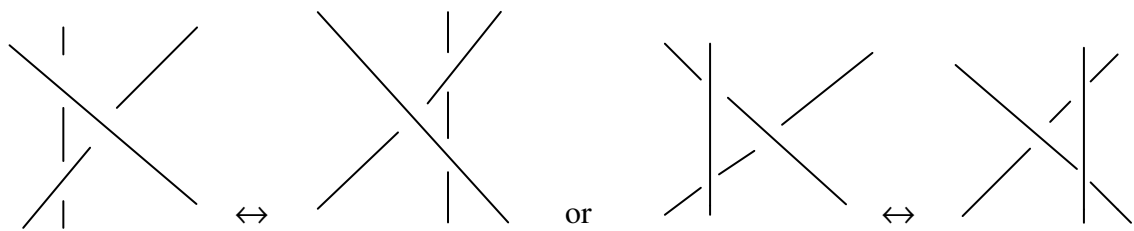

A more specific problem is The Unknotting Problem. The unknot is defined as the knot which has no crossings—it looks like a ring(see below). The problem is as follows: Given an arbitrary knot, can we determine whether or not it is equivalent to the unknot, under the three Reidemeister moves. Below are some examples of different forms of the unknot.

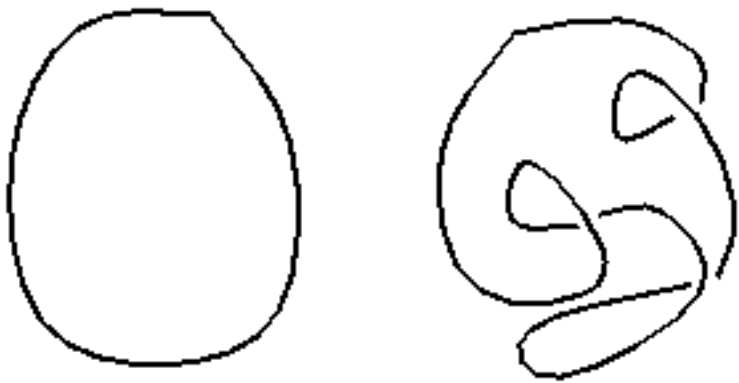

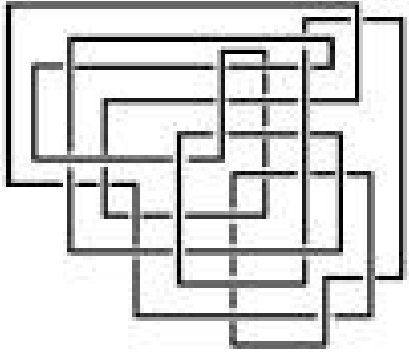

Note: Above picture copied from [J]

There are many different methods of approaching the Knot Equivalence Problem and the Unknotting Problem.

## 6.1.1 The Normal Surface Approach

In 1961, W. Haken[H] obtained a procedure to decide if a knot is unknotted using a theory of normal surfaces. He also outlined a similar approach to the Knot Equivalence Problem. In [HLP], this result of Haken was investigated and formalized into an algorithm to test for unknottedness. They proved that the Unknotting Problem is in NP. They also obtained a specific bound on the complexity for an algorithm to determine if a given knot with n crossings is unknotted—$O(\exp(cn))$.

In the above paper, the authors posed the problem of determining whether the Unknotting Problem is in co-NP. In other words, for a knot which is not equal to the unknot, can one produce a certificate of knottedness in polynomial time?

## 6.1.2 The Knot Group Approach

Another approach to solving these problems is by studying the knot group, a group which can be associated with any knot diagram. In [P], Papakyriakopoulos gave a rigorous proof of Dehn's

Lemma[D2] which states that if a knot group is infinite cyclic, then the associated knot is the unknot. Thus, studying the group of a given knot is a method of trying to solve the unknotting problem.

### 6.1.3 The Polynomial Invariant Approach

A different approach is to try to find polynomial invariants of knot classes. For instance, there is a polynomial called the Jones Polynomial which is defined and can be computed for all knots[J]. There is a conjecture that the Jones Polynomial of knot is equal 1 if and only if the knot is the unknot. This would, of course, solve the unknotting problem.

### 6.1.4 The Closed Braid Approach

Braid group theory provides us with yet another approach to these problems. In 1923, W. Alexander[A] proved that every knot can be represented as a closed braid. A closed braid is a braid which we "close up" by pulling the top bar of a braid around and attaching it to the bottom bar.

Example: Closure of 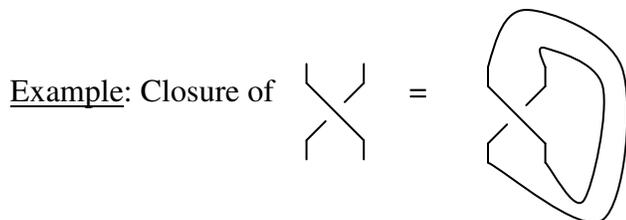

Thus instead of studying knots directly, we can instead study the corresponding braids. We can, thereby, study knots in a purely algebraic fashion. The Knot Equivalence Problem and the Unknotting Problem can now be reformulated as decision problems in group presentations. This approach is formalized in a result known as Markov's Theorem[M]. It states that for any two equivalent knots, the corresponding braids must be related by a chain of Markov moves. In other words, if β and γ are braids whose closures represent equivalent knots, then we can form a chain $\beta=\beta_1\rightarrow\beta_2\rightarrow\ldots\rightarrow\beta_i\rightarrow\ldots\rightarrow\beta_n=\gamma$, where each $\beta_{i+1}$ can be attained from $\beta_i$ from one of the following Markov moves:

1) If $\beta_i \in B_n$, then $\beta_{i+1}= \alpha^{-1}\beta_i\alpha$, for some $\alpha \in B_n$ (conjugation)

2) If $\beta_i \in B_n$, then $\beta_{i+1}= \beta_i\sigma_n^{\varepsilon}$, where $\varepsilon=1$ or $-1$ (stabilization)

or if $\beta_i=\beta\sigma_n^{\varepsilon}$, for $\varepsilon=1$ or $-1$ and $\beta \in B_n$, then $\beta_{i+1}=\beta$ (destabilization)

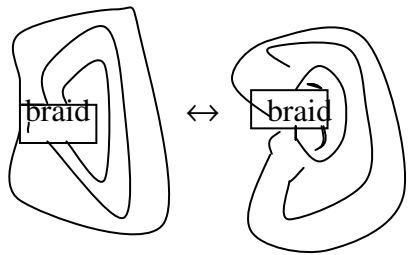

Markov Move 1 - conjugation

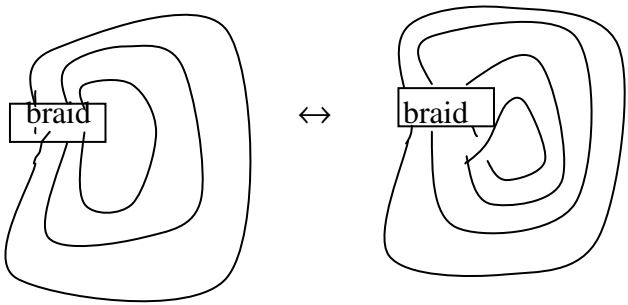

Markov Move 2 – stabilization/destabilization

Based upon this theorem, we are led to study the decision problems of section 1.3 in the case of the braid group. More specifically, the conjugacy problem acquires great importance. One must also consider exactly what effect an arbitrary sequence of conjugations and stabilizations/destabilizations have upon a braid word. Accurate answers to these questions can lead to algorithmic solutions to the Knot Equivalence Problem and to the Unknotting Problem.

As a furtherance of the above approach of using braids to study knots, Birman and Menasco wrote a series of papers studying the problem using braid foliations, a concept from geometric topology. This study culminated in the Markov Theorem Without Stabilization which was proven in [BM]. This theorem proves that in studying the braid representatives of two equivalent knots, we can form a Markov chain from one to the other without having to increase the braid index in the process. In this paper stabilization is replaced by other "complexity-reducing" moves.

Alongside of this approach, Birman and Hirsch [BH] developed an algorithm for solving the Unknotting Problem. They based this algorithm upon the fact that a knot is unknotted if and only if it is the boundary of a disc with a combinatorial foliation. No results have been proven as to the complexity of this algorithm.

## 6.2 Applications to Cryptography

One approach which is used in designing public-key cryptosystems is to utilize hard decision problems from combinatorial group theory. In 1999, [AAG] suggested such an approach by using groups whose word problem is easy, but whose conjugacy problem is very difficult. They suggested using the braid group for such a purpose. Another public-key cryptosystem using the braid group was introduced in [KLCHKP]. Thus, results about finding algorithmic solutions to decision problems in the braid group are sure to have applications in the field of cryptography. (see section 7 for more details)

## 6.3 P=NP?

As mentioned above, [PR] showed that the shortest word problem in the braid group is at least as hard as an NP-complete problem. Thus, if one could design an algorithm which would solve the shortest word problem in polynomial time, this would effectively prove that P=NP.

## 6.4 Using Computers for Computations in Braid Groups

There have been a number of computer programs that have been developed in order to facilitate the study of braid groups. Among these are: CBraid[C], BraidLink[Br], and Chevie[Ch].

# 7 Braid Group Cryptography

## 7.1 Some introductory remarks

One application of braid groups that is currently very active is that of braid group cryptography. The first paper that discussed the idea of using difficult problems from combinatorial group theory to design a public key cryptosystem is [WM]. However, the idea of using the braid group as a platform for cryptosystems was introduced in [AAG]. This paper introduced a protocol for a public-key cryptosystem involving braid groups. Later, in [KLCHKP], a different, Diffie-Hellman type of public-key cryptosystem was introduced, also using braid groups as its platform. Most recently, digital signature schemes which use braid groups as their platform have been introduced. [KCCL],[DGS]. As time goes on, the prospects of using braid groups in many areas of cryptography are slowly opening up.

The basic properties of braid groups which all of these applications make use of are (i) the fact that the word problem in the braid group can be quickly solved, (ii) elements in the braid group can quickly be put into a unique normal form, (iii) the conjugacy problem in the braid group offers no known, quick solution. In this section, we will study both of the above-mentioned cryptosystems. We will have some discussion of choices of parameters and then take a survey of the current attacks that have been launched against each of these cryptosystems. We will conclude with some directions for further research.

Definition: A *protocol* is a step-by-step method which is followed by two or more parties who are trying to achieve a certain objective

Definition: A *key-agreement protocol* is a protocol, which is used for the purpose of agreeing on a secret key, a shared secret, which will be used in a cryptosystem.

Definition: A *public-key-cryptosystem* (PKC) is an algorithmic method used for exchanging secure messages over an insecure channel. The assumption is that the parties involved do not share a key in advance. The key-agreement protocol is a central feature of any PKC.

The major PKC's in use today are based on finite abelian groups. They make use of hard problems in number theory. Examples of these are the Diffie-Hellman PKC, the RSA PKC, and the Elliptic Curve PKC. The braid groups, on the other hand, are infinite non-abelian groups. The two suggested public-key-cryptosystems using braid groups make use of hard problems in combinatorial group theory- specifically, the conjugacy problem in the braid group. We will illustrate below the two key-agreement protocols and will discuss the problems in braid group theory upon which they are based.

## 7.2  The Key Agreement Protocol of [AAG]

The objective of any key agreement protocol is for two parties, Alice and Bob, to attain a shared secret key, which cannot be computed by an outsider who views all transmissions between Alice and Bob. The shared key in the key-agreement protocol of [AAG] is the commutator of a and b, two elements of $B_n$. Namely, $[a,b] = aba^{-1}b^{-1}$. We begin by making public (i) the braid index $n \in \mathbf{N}$. (ii) a subgroup of $B_n$, say $G_A=<a_1,\ldots,a_r>$, and (iii) a subgroup of $B_n$, say $G_B=<b_1,\ldots,b_s>$.

Now, we will trace Alice's actions:

i) Alice chooses a word a in the subgroup $G_A$. She then conjugates each $b_i$ by a to attain $ab_1a^{-1},\ldots,ab_sa^{-1}$. Call these elements $b_1',\ldots,b_s'$, respectively.

ii) Alice then computes the normal form for each of the elements $b_i'$ and transmits them to Bob.

We will now trace Bob's actions, which are similar to Alice's:

i) Bob chooses a word b in the subgroup $G_B$. He then conjugates each $a_i$ by b to attain $ba_1b^{-1},\ldots,ba_rb^{-1}$. Call these elements $a_1',\ldots,a_r'$, respectively.

ii) Bob then computes the normal form for each of the elements $a_i'$ and transmits them to Alice.

Using the information that they have, both Alice and Bob can now compute the commutator of a and b. This is possible because of the fact that conjugation is a homomorphism. Namely, the product of the conjugates of two elements is equal to the conjugate of the product. i.e. $(zxz^{-1})(zyz^{-1})=z(xy)z^{-1}$. This being the case, assume Alice chose $a=a_{i1}a_{i2}\ldots a_{ik}$, the product of k elements in $G_A$. Then $bab^{-1}=ba_{i1}a_{i2}\ldots a_{ik}b^{-1}=(ba_{i1}b^{-1})(ba_{i2}b^{-1})\ldots(ba_{ik}b^{-1})= b_{i1}'b_{i2}'\ldots b_{ik}'$, which Alice can compute using the values of $b_i'$ that Bob sent to her and the $a_{i\_}$'s which she chose. Thus Alice computes:

$[a,b]= aba^{-1}b^{-1} = a(bab^{-1})^{-1}$.

In a similar manner, Bob can compute $aba^{-1}$ using the values of $b_i'$ sent to him by Alice. He thereby computes $[a,b]= aba^{-1}b^{-1}$. Thus, both Alice and Bob have a representation for the commutator of a and b. If they each convert their result to normal form, then they have a shared secret key. This is the key-agreement protocol of [AAG]. Notice that the strength of this protocol lies in the difficulty of solving the multiple simultaneous conjugacy problem. We will discuss this further when we discuss attacks.

In [AAFG], the authors improve on this protocol by introducing an efficient key extractor, using the Colored Burau Group.

## 7.3 The Key-Agreement Protocol of [KLCHKP]

This scheme is the braid group equivalent of the Diffie-Hellman key agreement protocol. The shared key in this protocol is $(ab)^{-1}x(ab)$, as will be explained. First Alice and Bob make public (i) a braid index $n \in \mathbf{N}$, (ii) integers $p,r \in \mathbf{N}$ such that l+r=n. (iii) a sufficiently complicated braid word $x \in B_n$. Now, we denote by $LB_p$ the subgroup of $B_n$ generated by $\sigma_1,\ldots,\sigma_{p-1}$, and by $RB_r$ the subgroup of $B_n$ generated by $\sigma_{n-r+1},\ldots,\sigma_{n-1}$. Now Alice and Bob follow the following steps:

i) Alice chooses a secret element $a \in LB_l$. She then computes $axa^{-1}$, converts it to normal form and sends it to Bob.

ii) Bob chooses a secret element b∈ $RB_r$. He then computes $bxb^{-1}$, converts it to normal form and sends it to Alice.

iii) Alice takes the $bxb^{-1}$ that was sent by Bob and computes $abxb^{-1}a^{-1}$.

iv) Bob takes the $axa^{-1}$ that was sent by Alice and computes $baxa^{-1}b^{-1}$.

Now by the defining relations in a braid group, all elements from $LB_l$ commute with all elements in $RB_r$. Thus, a and b commute. Therefore, Alice and Bob have attained equivalent elements. If they both convert their result to normal form, then they'll have a shared secret key. Notice that the strength of this scheme lies in the difficulty in solving the generalized conjugacy search problem. (this is the conjugacy problem, with the restriction that any conjugating element can only come from a specified subgroup of the initial group, in this case $LB_p$ or $RB_r$ ). We will discuss this further when we consider attacks. Using this key agreement protocol, the authors proceeded to develop a new PKC.

In [CKLHC], a revision of this scheme was proposed, in which this scheme was made more general. Additionally, specific algorithms were designed to implement this cryptosystem

## 7.4 Discussion of Parameters and Attacks on [AAG]

In order to implement any cryptosystem, the parameters must be chosen in a manner which will ensure security against all attacks. The parameter which were suggested in [AAFG] are as follows: (i) choose braid index n=80 or larger. (ii) Let the public subgroups $G_A$ and $G_B$ be generated by 20 elements each, where each generator is the product of 5 to 10 Artin generators. Also, let each set of generators make use of all the Artin generators in $B_n$ (iii) Let the secret words a and b each have length at least 100 in the public generators.

Since this cryptosystem has been developed, there have been a number of attacks that must be dealt with. We will mention a few of them now.

i)         In [H], Hughes used the Burau representation of the braid group to mount a linear algebraic attack for certain classes of keys used in this cryptosystem. His conclusion was that in order for this scheme to be successful, one must be more specific in the choice of keys.

ii)        In [HT], Hughes and Tannenbaum described a probabilistic "length" attack on this PKC. This attack makes use of a length function on the set of conjugates defining the public key. Their conclusion was that to avoid this attack, the generators of the subgroups chosen must have small canonical length. This, they claim, may make the system more vulnerable to an attack similar to (i).

iii)       In [LL], the authors put forth a linear algebraic attack on the key extractor used in [AAFG]. They also describe a mathematical algorithm to solve the multiple simultaneous conjugacy problem in braid groups. Such a solution would attack the private key used in this cryptosystem. They suggest that by selecting generators of the public subgroups to be pure braids, this attack may be avoided.

iv)       In [Go], Gonzalez-Meneses improved on [LL]'s algorithm to solve the multiple simultaneous conjugacy problem in braid groups.

v)        In [HS], the authors describe an algorithm for solving certain instances of the conjugacy problem in the braid group. They claim that the modifications which have been suggested to avoid the attacks of [HT] and [LL], will not provide security against their attack.

### 7.5 Discussion of Parameters and Attacks for [KLCHKP]

The parameters suggested for the [KLCHKP] are as follows: (i) Choose braid index, n=45. (ii) the public x should be composed of about 1450 Artin generators and should involve all the Artin generators in $B_n$ (iii) the private a and b should be composed of about 360 Artin generators each. (these parameters are such that all of the braids chosen will have canonical length at least 3). The following attacks have been launched against this PKC:

i) In [H2], Hughes utilized the left super-summit-set, an invariant under conjugation, to attack this PKC with the given parameter choices.

ii) In [HS], the authors attack this system for even larger parameter choices. They use an algorithm which they designed for the solution to certain instances of the conjugacy problem in braid groups. However, this attack cannot be launched against the generalized version of this PKC, which was proposed in [CKLHC]. In this revised version, the underlying problem comes down to the following: Given $a$ and $x_1 a x_2$, find any $z_1$ and $z_2$ such that $z_1 a z_2 = x_1 a x_2$, where $a \in B_n$ and $x_1, x_2, z_1, z_2 \in LB_p$.

iii) In [LP], the authors suggest a method of solving the above mentioned problem in the generalized cryptosystem. They followed Hughes[H1] lead in using the Burau representation to solve the conjugacy problem in the braid group. (See Section 1.3.3 for a basic discussion of the linear algebraic method of solving the conjugacy problem in braid groups). They proposed two algorithms, one which improves on Hughes' in efficiency, the other in accuracy. In order to break the PKC of [CKLHC], it is not necessary to solve the conjugacy problem, but only the following problem: Given $a, b \in B_n^+$, find $x_1, x_2 \in LB_n^+$ satisfying $b = x_1 a x_2$, provided that such braids exist. The authors solve this problem by (i) transforming it into a problem in $GL_k(R)$, using the Burau representation, (ii) Solving it in $GL_k(R)$, and then (iii) lifting the solution back into $B_n$. Although the Burau representation is not faithful, they still utilize it as opposed to the Lawrence-Krammer representation which is faithful. The reason for this is because under the Lawrence-Krammer representation, a braid in $B_n$ is sent to a matrix of size $n(n-1)/2 \times n(n-1)/2$. Such large matrices seem too large to be practical to work with. In contrast, under the Burau representation, a braid in $B_n$ is sent to an n x n matrix. Although this is more manageable in size, there seems to be one major difficulty in using the Burau representation. Namely, it is not faithful. Thus, when why try to lift a solution from $GL_n(R)$ back to $B_n$, we may not get the correct braid. Despite this fact, Hughes designed an algorithm which was successful in a similar task with a very high probability. In order to break a PKC, we do not have to be correct 100% of the time--we simply need to make the PKC insecure. Therefore, the authors analyzed Hughes' idea and its success and, thereby, improved upon it. Using their solution, the authors attempted to break the PKC using several parameter choices suggested by the inventors of the PKC. Their results were as follows: For seven of nine parameter choices, this method successfully recovered the private key from the public key with significant probability and in a reasonable time. They also gave a requirement for secure parameters against their attack which, as it turns out, conflict with the parameters which will be secure against a brute-force attack. Despite the force of this attack, the authors maintain that it is inefficient for very large parameters as those suggested in [CKLHC]. They, therefore, conclude that at the current state of knowledge, the revised cryptosystem is still secure against all attacks.

iv) In [CJ], the authors propose another linear algebraic attack on the PKC of [KLCHKP]. However, this approach does not use the Burau representation, but uses the Lawrence-Krammer representation for $B_n$, which is faithful. They use this to solve the Diffie-Hellman conjugacy problem (a case of the generalized conjugacy problem) and, thereby, attack the original version of the PKC in [KLCHKP]. Although the image of the Lawrence-Krammer representation consists of matrices which are too large to deal with in solving the conjugacy problem, the authors claim that when we restrict our attention to the Diffie-Hellman conjugacy problem, the solution becomes more manageable. They also show that in order to break the encryption scheme and the key agreement scheme, it is not necessary to recover the original key. It suffices to find a "fake" key which plays the same role in these schemes as the real key does. Additionally, the authors claim that the Lawrence-Krammer representation is invertible in polynomial time in the braid index and the length of the word involved. Using all of these facts, together with linear algebraic calculations, the authors conclude with the following theorem.

<u>Theorem</u>: Assume $LB_n$ and $RB_n$ are two commuting subgroups of the 2n-braid group $B_{2n}$. Given $u \in B_n$, $a^{-1}ua$, $b^{-1}ub$, for $a \in LB_n$ and $b \in RB_n$, $b^{-1}a^{-1}uab$ can be computed in $O(n^{45} \delta^{10} \log^2 n)$ bit operations where $\delta$ is the maximum word length of u, $a^{-1}ua$, $b^{-1}ub$.

The authors explain that although the complexity of this algorithm is polynomially bounded, it is too large to implement for the proposed parameters. However, they expect that the complexity can be reduced by a more precise analysis. It seems that there is along way to go to make this into a practical attack. Another limitation of this attack is that it seems to be limited to the original version of the PKC from [KLCHKP], as was the attack of [HS]. However, it seems unable to attack the revised version of [CKLHC].

It would be interesting to apply the approach of [CJ] to other algorithmic problems in the braid group. (See Open Problem 7).

## 7.6 Authentication and Signature Schemes Using Braid Groups

The latest application of braid groups to the area of cryptography is in the creation of signature or authentication schemes based upon braids. The first such suggestion is found in [KCCL]. The authors devised a digital signature scheme which is based upon the fact that in the braid group, the conjugacy search problem is hard while the conjugacy decision problem is feasible. Then in [DGS], three public-key authentication schemes based upon the braid group were presented. The first of these schemes is based upon the Diffie-Helman type of key agreement protocol which was developed in [KLCHKP]. The second and third schemes are what are known as *zero knowledge interactive proofs of knowledge.* The second scheme is based upon the difficulty of the conjugacy search problem in braid groups. The third scheme, in addition to using the difficulty of the conjugacy search problem, also makes use of the difficulty of the root problem in the braid groups. Namely, assuming that a braid x is a $p^{th}$ power in $B_n$, it is difficult to find a braid $y \in B_n$ such that $x=y^p$, i.e it is difficult to find the $p^{th}$ root of x in $B_n$. In this paper, it was also suggested that instead of using normal forms to implement these schemes, it may be more efficient to use the method of the handle reduction algorithm described in [De].

## 7.7  Concluding Remarks

With the passing of every month, it seems that there are two types of papers published regarding the use of braid groups as a platform for modern cryptosystems. On the one hand, there are many new and interesting methods being developed to apply braid groups to different cryptographic problems. At the same time, one finds an equal number of researchers pointing out weaknesses in such an approach, and attacking these cryptosystems. It seems that this is the fate of any area of cryptography. The pursuit of attaining security in this insecure world is a topic bound to draw the interest of laymen and the intelligence of mathematicians. It would seem that the only way to unravel the truth in this area is to let time run its course. I am glad to be involved in the area of braid group algorithmics and to be able to be a part of this volatile world of braid group cryptography in a first hand manner. I would like to pursue further research in attempts to validate the usage of braid groups in modern cryptography.

# 8 Open Questions

1) As mentioned above, in [HLP] it is proven that the unknotting problem is in NP. In other words, for a given representation of the unknot, one can produce a certificate of unknottedness in polynomial time. However, the proof given is of a geometric nature, and involves a complicated construction. It would be desirable to find a more concrete, algebraic method to prove this result. It seems that such a proof should exist. I have worked on this problem, but have not found a proof.(see open problem 8)

2) A similar problem, also raised in [HLP] is to determine if the unknotting problem is in co-NP. In other words given a nontrivial knot, can one produce a certificate of knottedness in polynomial time? This is an open problem. This problem is interesting because there are but a few problems known to be in NP $\cap$ co-NP and not known to be in P. (see open problem 8)

3) In [St], an algorithm is presented to find if a given braid $\beta \in B_n$ has an $m^{th}$ root, for some $m \in \mathbf{N}$. In other words, one can decide if there exists an $\alpha \in B_n$ such that $\alpha^m = \beta$.(see Section 2.8)  It would be interesting to improve upon the results in this paper. Firstly, can we find a bound on the complexity of the above algorithm? ( all that is proven in [St] is that the problem is solvable) Additionally, for a given $\beta$ and a given $m \in \mathbf{N}$, can we find every $m^{th}$ root of $\beta$ of bounded length?

4) In section 2, we found a quick solution to the generalized word problem in the case of cyclic subgroups of the braid group where exp of the generator is nonzero. Can we find an analogous solution to the case of cyclic subgroups of the braid group which are generated by an element whose exp is zero?

5) Can we find a bound on the complexity of the conjugacy problem in the cyclic amalgamation of two braid groups? Above, we have proven that this problem is solvable, but it would be nice to have a good bound. It would seem that this problem comes down to finding a bound on the solution to the following problem: "Given $\alpha, \beta \in B_n$, and given $\sigma_i \in B_n$, determine if there exists a $c \in \mathbf{Z}$ such that $\sigma_i^c \alpha \sigma_i^{-c} = \beta$. If there exists such a c, find it." We showed in Section 5.3 that this is solvable, but we have not yet found a good bound. This problem seems quite simple, but is elusive.

6) In this paper we considered the construction of the free product of braid groups amalgamating a cyclic subgroup. We considered decision problems in this setting. What if we were to change the construction? For instance, let's say instead of amalgamating a cyclic subgroup, we amalgamated some other subgroup, i.e., a finitely generated subgroup. Then what would happen with the word and conjugacy problems? Or, let's say we considered an HNN-extension of the braid group. Then what could we say regarding these decision problems? This would be a nice direction for future research.

7) Can we make any improvement upon the bound for the solution to the conjugacy problem in the braid group? Can we verify BKL's conjecture that it is in fact in polynomial time? If we cannot make any advances in this problem, perhaps we can improve our results regarding related problems. Namely, if we can find a quick algorithm to solve the "Diffie-Hellman conjugacy problem" in braid groups, then we could break the PKC of [KLCHKP], and other Diffie-Helman-like applications of the braid group to cryptography. Similarly, if we can find a quick solution to the multiple simultaneous conjugacy problem, we could break the PKC of [AAG]. Thus, any improvements in our knowledge regarding different variants of the conjugacy problem in the braid group will certainly be significant from a theoretical and a practical standpoint.

8) In [CJ], the authors utilized the Lawrence-Krammer representation of the braid group to give a solution to the Diffie-Hellman conjugacy problem in braid groups. They claimed that this representation is invertible in polynomial time. Thereby, they solved the problem in the image of this faithful representation, and then lifted it back to the braid group. This opens up a general approach to solving algorithmic problems in the braid group. We can always convert any such problem into the analogous problem in the image of the representation. It would be a nice research project to go through different open problems in the braid group and apply such a method. Perhaps, we would be able to solve formerly unsolved problems, or at least improve upon certain known complexity bounds. With specific reference to results in our paper, I would like to apply this method to solving open problem (5) and, thereby, be able to find a good complexity bound for the conjugacy problem in the cyclic amalgamation of two braid groups. Additionally, it would be interesting to study the Knot Equivalence Problem and the Unknotting Problem using such a method. The approach would be to consider what happens to a sequence of Markov moves which are applied to a given braid, when considered in the image of the Lawrence-Krammer representation. Perhaps we could find another proof that the Unknotting Problem is in NP, or prove that it is in co-NP.(see open problems 1 and 2)

9) Can we improve upon the known applications of braid groups to cryptography? Can we find new applications? Or, perhaps, can we show that some or even all of the suggested applications of braid groups to cryptography do not stand up to the demands of security which are necessary in today's world?